\documentclass[3p]{elsarticle}

\usepackage[all]{xy}
\usepackage{bibgerm}
\usepackage[english]{babel}

\usepackage[active]{srcltx}
\usepackage{amsmath}
\usepackage{amsthm}
\usepackage{amsfonts}
\usepackage{amssymb}

\newtheorem{defi}{Definition}[section]
\newtheorem{remark}[defi]{Remark}

\newtheorem{thm}[defi]{Theorem}
\newtheorem{lemma}[defi]{Lemma}
\newtheorem{corollary}[defi]{Corollary}
\newtheorem{notation}[defi]{Notation}
\newtheorem{prop}[defi]{Proposition}

\newcommand{\Q}{\mathbb{Q}}
\newcommand{\Z}{\mathbb{Z}}
\newcommand{\N}{\mathbb{N}}
\newcommand{\F}{\mathbb{F}}
\newcommand{\SL}{\operatorname{SL}}
\newcommand{\Sl}{\operatorname{SL}}
\newcommand{\OO}{\mathcal{O}}
\newcommand{\modC}{\operatorname{\mathbf{mod}}}
\newcommand{\projC}{\operatorname{\mathbf{proj}}}

\newcommand{\iso}{\cong}
\newcommand{\trred}{\operatorname{tr. red.}}
\newcommand{\eps}{\varepsilon}
\newcommand{\Lifts}{\ensuremath{\widehat{\mathfrak{L}}}}
\newcommand{\Out}{\operatorname{Out}}
\newcommand{\Aut}{\operatorname{Aut}}
\newcommand{\id}{{\operatorname{id}}}
\newcommand{\Jac}{\operatorname{Jac}}
\newcommand{\Gal}{\operatorname{Gal}}
\newcommand{\Rad}{\operatorname{Rad}}
\newcommand{\End}{\operatorname{End}}
\newcommand{\Hom}{\operatorname{Hom}}
\newcommand{\Soc}{\operatorname{Soc}}
\newcommand{\length}{\operatorname{length}}
\newcommand{\diag}{\operatorname{diag}}
\newcommand{\rank}{\operatorname{rank}}
\newcommand{\op}{\operatorname{op}}
\newcommand{\Res}{\operatorname{Res}}
\newcommand{\add}{\operatorname{add}}
\newcommand{\Ext}{\operatorname{Ext}}

\title{The $p$-adic group ring of $\SL_2(p^f)$}
\date{}

\begin{document}

\begin{frontmatter}
 \author{Florian Eisele}
\ead{feisele@vub.ac.be}
\address{Vrije Universiteit Brussel, Pleinlaan 2, 1050 Brussel, Belgium}
\begin{abstract}
    In this article we show that the $\Z_p[\zeta_{p^f-1}]$-order $\Z_p[\zeta_{p^f-1}]\SL_2(p^f)$ can be recognized among those orders whose reduction modulo $p$ is isomorphic to $\F_{p^f}\SL_2(p^f)$ using only ring-theoretic properties (in other words we show that  $\F_{p^f}\SL_2(p^f)$  lifts uniquely to a  $\Z_p[\zeta_{p^f-1}]$-order, provided certain reasonable conditions are imposed on the lift).
    This proves a conjecture made by Nebe in \cite{NebeSl2Char2} and \cite{NebeSl2Charp}
    concerning the basic order of $\Z_p[\zeta_{p^f-1}]\SL_2(p^f)$.
\end{abstract}
\end{frontmatter}

\section{Introduction}

Let $p$ be a prime and let $(K,\OO,k)$ be a $p$-modular system. This article is concerned with the group ring
$\OO \SL_2(p^f)$ for some $f\in \N$. Hence we are dealing with the discrete valuation ring version of what is typically referred to as representation theory in ``defining characteristic''.
Our aim in this paper is to prove a conjecture made by Nebe in \cite{NebeSl2Char2} (for the case $p=2$) respectively 
\cite{NebeSl2Charp} (for the case $p$ odd) which claims (rightly) to describe the group ring of $\SL_2(p^f)$ over sufficiently large extensions $\OO$ of $\Z_p$. Here, ``to describe the group ring'' means to describe its basic order.
However, our proof of Nebe's conjecture is indirect, and consists essentially of showing that a ``unique lifting theorem'' (see Corollary \ref{corollary_unique_lifting}) holds for the group ring of $\SL_2(p^f)$. Basically this unique lifting theorem asserts that (provided $k \supseteq \F_{p^f}$) any $\OO$-order reducing to $k \SL_2(p^f)$ which has semisimple $K$-span and is self-dual (with some technical condition on the bilinear form with respect to which it is self-dual) has to be isomorphic to $\OO \SL_2(p^f)$.
Nebe's conjecture is an immediate consequence of this, but the theorem may well be considered an interesting result in its own right.

This work is a continuation of the author's work in \cite{EiseleDerEq}, where a ``unique lifting theorem'' similar to the one mentioned above is proved for $2$-blocks with dihedral defect group. Our approach is, as in \cite{EiseleDerEq}, based on the idea that, provided it is properly formulated, such a theorem holds for a $k$-algebra  if and only if it holds for all $k$-algebras derived equivalent to the original one. By the abelian defect group conjecture (which is known to be true in the special case encountered in the present paper), the blocks of $k\SL_2(p^f)$  are derived equivalent to their Brauer correspondents (we must assume $k$ to be algebraically closed for this, but we manage to work around that). And, as it turns out, to prove a ``unique lifting theorem'' for these Brauer correspondents is fairly easy due to their simple structure. In particular we prove Nebe's conjecture without ever having to put up with the complicated combinatorics that arises in the representation theory of 
$\SL_2(p^f)$.

\section{Notation and technical prerequisites}

Throughout this article, $p$ will denote a prime and $(K,\OO, k)$ will denote a $p$-modular system such that $K$ is a complete and unramified 
extension of $\Q_p$. We let $\bar K$ and  $\bar k$ denote the respective algebraic closures. By $\nu_p:\ K \longrightarrow \Z$ we denote 
the $p$-valuation on $K$.
\begin{notation}
 We are going to use the following notations (all of which are more or less standard):
\begin{itemize}
\item $\modC_A$ and $\projC_A$: the categories of finitely-generated modules respectively finitely-generated projective modules over the ring $A$.
\item $\mathcal D^b(A)$, $\mathcal D^-(A)$: the bounded respectively right bounded derived category of $A$-modules.
\item $\mathcal K^b(\projC_A)$: the homotopy category of bounded complexes with finitely generated projective terms. 
\item $-\otimes^{\mathbb L}_A =$: the left derived tensor product.
\item $\Out_k(A)$: the outer automorphism group of the $k$-algebra $A$. To keep notation simple we will not differentiate between elements of 
 $\Out_k(A)$ and representatives for those elements in $\Aut_k(A)$.
\item $\Out_k^0(A)$ (assuming $k$ is algebraically closed): the identity component of the algebraic group $\Out_k(A)$.
\item $\Aut_k^s(A)$ and $\Out_k^s(A)$: These denote the subgroups of $\Aut_k(A)$ respectively $\Out_k(A)$ which stabilize all isomorphism 
classes of simple $A$-modules (with the action of $\Aut_k(A)$ and $\Out_k(A)$ on isomorphism classes of modules being given by twisting). 
\item If $A$, $B$ and $C$ are rings, and $\alpha:\ A\rightarrow C$ as well as $\beta:\ B\rightarrow C$ are ring homomorphisms, then 
we denote by ${_\alpha}C_\beta$ the $A$-$B$-bimodule which as a set coincides with $C$, where $a\in A$ and $b\in B$ act 
on $c\in C$ by the formula $a\cdot c \cdot b := \alpha(a)\cdot c \cdot \beta(b)$.
\end{itemize}
\end{notation}

We set $\SL_2(p^f) := \SL_2(\F_{p^f})$ and
\begin{equation}
\Delta_2(p^f) := \left\{  \left[\begin{array}{cc} a&b\\0&a^{-1}\end{array}\right]\ \bigg|\ a,b\in \F_{p^f},\ a\neq 0 \right\} \iso C_p^f \rtimes C_{p^f-1}
\end{equation}
Note that $\Delta_2(p^f)$ is the normalizer of a $p$-Sylow subgroup of $\SL_2(p^f)$, namely of the group of unipotent
upper triangular $2\times 2$-matrices. Also note that $k$ splits $\SL_2(p^f)$ and $\Delta_2(p^f)$ if and only if $k\supseteq \F_{p^f}$.

One important property of group rings over integral domains which we are going to exploit in this article is that they are self-dual with respect to a bilinear form of the kind defined in the following definition.
\begin{defi}[Trace bilinear form]
	Let 
	\begin{equation}
		A = \bigoplus_{i=1}^l D_i^{n_i\times n_i}
	\end{equation}
	be a finite-dimensional semisimple $K$-algebra given in its Wedderburn decomposition (i. e. the $D_i$ are division algebras over $K$ and the $n_i$ are certain natural numbers). Given an element $u=(u_1,\ldots,u_l) \in Z(A)=Z(D_1)\oplus\ldots\oplus Z(D_l)$ we define a map
	\begin{equation}\label{eqn_ddjk33jkljk}
		T_u:\ A \longrightarrow K: \ a=(a_1,\ldots,a_l) \mapsto \sum_{i=1}^l {\rm tr}_{Z(D_i) / K}\trred_{D_i^{n_i\times n_i}/Z(D_i)} (u_i\cdot  a_i)
	\end{equation}
	and (by abuse of notation) a bilinear form of the same name: $T_u:\ A\times A \longrightarrow K:\ (a,b)\mapsto T_u(a\cdot b)$.
	Here ``$\rm tr_{Z(D_i)/K}$'' denotes the trace map in the sense of Galois theory, and 
	``$\trred_{D_i^{n_i\times n_i}/Z(D_i)}$'' denotes the reduced trace as defined for central simple algebras. 
	
	For a full $\OO$-lattice $L\subset A$ we define its \emph{dual} as follows
	\begin{equation}
		L^{\sharp,u} := \{ a \in A \mid T_u(a, L) \subseteq \OO \}
	\end{equation}
	We call $L$ self-dual (with respect to $T_u$) if $L^{\sharp,u} = L$ (the ``$u$'' may be omitted when its choice is clear from context).
\end{defi}

\begin{remark}\label{remark_symm_elem}
	\begin{enumerate}
		\item The definition of $T_u$ as given above is compatible with extensions of scalars in the following sense:
		If $K'$ is a field extension of $K$, $\OO'$ is the integral closure of $\OO$ in $K'$ and $\Lambda$ is a full $\OO$-order in the semisimple $K$-algebra $A$, then 
		$\Lambda$ is self-dual in $A$ with respect to $T_u$ if and only if $\OO'\otimes\Lambda$ is self-dual in $K'\otimes A$ with respect to 
		$1\otimes u$. Therefore we will often think of $u$ as an element of $Z(\bar K \otimes A)$.
		\item An order $\Lambda \subset A$ is self-dual with respect to some form $T_u$ if and only if $\Lambda$ is a symmetric $\OO$-order (but of course,
		the element $u\in Z(A)$ such that $\Lambda = \Lambda^{\sharp,u}$ contains more information than merely that the order in question is symmetric). 
		\item Group rings $\OO G$ (for finite groups $G$) are self-dual orders. Let $\chi_1,\ldots,\chi_l$ denote the (absolutely) irreducible $\bar K$-valued 
		characters of $G$. Hence
		\begin{equation}
		\bar K G \iso 	\bigoplus_{i=1}^l \bar K^{\chi_i(1)\times \chi_i(1)}
		\end{equation}
		is the Wedderburn decomposition of $\bar K G$. Then  $\OO G = \OO G^{\sharp,u}$, where
		\begin{equation}
			u = \left ( \frac{\chi_1(1)}{|G|},\ldots, \frac{\chi_l(1)}{|G|}\right ) \in Z(KG) \subset Z(\bar K G) \iso \bigoplus_{i=1}^l \bar K
		\end{equation}
	\end{enumerate}
\end{remark}

We will be using the following definition of decomposition numbers:
\begin{defi}
    Let $\Lambda$ be an $\OO$-order with semisimple $K$-span.
    The decomposition matrix of $\Lambda$ is a matrix whose rows are labeled by the isomorphism classes of 
    simple $K\otimes \Lambda$-modules and whose columns are labeled by the isomorphism classes of simple $\Lambda$-modules.
    If $S$ is a simple $\Lambda$-module, $P$ is the projective indecomposable $\Lambda$-module with top $S$ and $V$ is a simple $K\otimes\Lambda$-module, then we define
    the entry $D_{V,S}$ to be the multiplicity of $V$ as a direct summand of $K\otimes P$.
\end{defi}

\section{Koshita's and Nebe's descriptions of the group ring}

In this section we are going to have a quick look at the descriptions of the basic algebra of the group algebra of $\SL_2(p^f)$
as given by Koshita and later, in the $p$-adic case, by Nebe. Our aim is to explain how to write down explicitly the description 
of the basic order of $\OO \SL_2(p^f)$ conjectured in \cite{NebeSl2Char2} (assuming as known the combinatorial description of the decomposition matrix of this order given in \cite{DecompSL2}), and to exhibit exactly which parts of it were actually of conjectural nature.
This is technically not a prerequisite to 
understanding the rest of this paper, since we will be dealing exclusively with the Brauer 
correspondents of the blocks of $k \SL_2(p^f)$. For simplicity's sake we will restrict our attention to the case $p=2$ (the case of odd $p$ works similarly, but happens to be a bit lengthier). 

In \cite{KoshitaChar2} respectively \cite{KoshitaCharp}, Koshita gave a description of 
the basic algebra of $\bar k \SL_2(p^f)$ as quiver algebra modulo relations, using the description of the projective 
indecomposable $\SL_2(p^f)$-modules given in \cite{AlperinProjSL2} as his starting point. Koshita's presentation is given in Theorem \ref{thm_desc_koshita} below.
\begin{notation} Let $N$ be a set and let $X,Y \subseteq N$ be subsets. Then we let $X+Y$ denote the symmetric difference between $X$ and $Y$
(i. e., $X+Y = X\cup Y - X \cap Y$).
\end{notation}
\begin{thm}[Koshita]\label{thm_desc_koshita}
	Let $Q$ be the quiver defined as follows:
	\begin{enumerate}
	\item the vertices of $Q$ are labeled by the subsets of $N := \Z/f\Z$.
	\item for any $I \subseteq N$ and any $i \in N$ such that $i-1 \notin I$ there is an arrow
		$\alpha_{i,I}: I + \{ i \} \longrightarrow I$.
	\end{enumerate}
	
	Then the basic algebra of $\bar k \SL_2(2^f)$ is isomorphic to the quotient of $\bar kQ$ by the ideal generated by the following 
	families of elements:
	\begin{enumerate}
	\item $\alpha_{i,I} \cdot \alpha_{j, I+\{i\}}-\alpha_{j,I} \cdot \alpha_{i, I+\{j\}}$ where $i-1$ and $j-1$ are not in $I$ and
	$j \notin \{i-1,i,i+1\}$
	\item $\alpha_{i,I}\cdot \alpha_{i,I+\{i\}}$ where $i$ and $i-1$ are not in $I$.
	\item $\alpha_{i+1,I}\cdot \alpha_{i,I+\{i+1\}} \cdot \alpha_{i,I+\{i\}+\{i+1\}}-\alpha_{i,I}\cdot \alpha_{i,I+\{i\}} \cdot \alpha_{i+1,I}$ where $i-1$ and $i$ are not in $I$.
	\item $\alpha_{i, I+\{i+1\}} \cdot \alpha_{i+1, I+\{i, i+1\}}\cdot \alpha_{i,I+\{i\}}$
	where $i\in I$ but $i-1\notin I$.
	\end{enumerate}
\end{thm}

\begin{defi}
	We denote the $\bar k$-algebra constructed in the foregoing theorem by $\bar \Lambda$. Moreover we let
	$\{\bar e_I\}_{I\subseteq N}$ be a system of pair-wise orthogonal primitive idempotents (where the indices correspond to the  respective vertices in $Q$ that the idempotents are associated with). For $I,J \subseteq N$ we define
	$\bar \Lambda_{IJ} := \bar e_I \bar \Lambda \bar e_J$. 
\end{defi}

\begin{remark}
	While our notation for the arrow $\alpha_{i,I}$ specifies the vertex from which it originates, this information is
	usually redundant when specifying a path, since the origin of an arrow must coincide with the target of 
	the arrow preceding it in the path. Therefore we make the following notational convention:
	\begin{equation}
		\alpha_i := \sum_{I \subseteq N-\{i-1\}} \alpha_{i,I}
	\end{equation}
\end{remark}
\newcommand{\weird}{\mathcal}
In \cite{NebeSl2Char2}, Nebe describes an $\OO$-order which reduces to a $k$-algebra with quiver and relations as in the foregoing theorem. The constructed order is self-dual, and its $K$-span is semisimple.
We will now outline this description. We assume for the remainder of this section that 
$\OO$ is an (unramified) extension of $\Z_2[\zeta_{2^f-1}]$, in order to ensure that both $k$ and $K$ are splitting fields for the group  $\SL_2(2^f)$.

Let $\weird R$ be the power set of $N = \Z / f \Z$.
As seen in Theorem \ref{thm_desc_koshita} the elements of $\weird R$ are in bijection with the (isomorphism classes of) simple $\bar k \SL_2(2^f)$-modules.
 Let $\weird C$ be an index set 
for the irreducible ordinary representations of $\SL_2(2^f)$. For $R \in \weird R$ denote by $\weird  C_R$ the subset of $\weird C$ corresponding to the irreducible ordinary representations which have non-zero decomposition number with the simple module associated with $R$. In the same vein, 
given $C\in \weird C$, define $\weird R_C$ to be subset of $\weird R$ corresponding to simple modules having non-zero decomposition number with the irreducible ordinary representation associated with $C$. Then the basic order of $\OO \SL_2(2^f)$ -- which we henceforth will refer to as $\Lambda$ -- is a full $\OO$-order in the split semisimple $K$-algebra
\begin{equation}
	A := \bigoplus_{C \in \weird C} K^{\weird R_C\times \weird R_C}
\end{equation}
We may assume that we have a complete set $\{e_R\}_{R\in \weird R}$ of pair-wise orthogonal primitive idempotents in $\Lambda \subset A$ such that each $e_R$
is diagonal in each of the matrix rings $K^{\weird R_C\times \weird R_C}$. The fact that all decomposition numbers of $\SL_2(2^f)$ are either zero or one implies that $e_R$ is simply a diagonal matrix unit in the direct summands of $A$ labeled by the elements of $\weird C_R$. Consequently, 
$\Lambda_{RR} := e_R \Lambda e_R$ is a commutative $\OO$-order, whose $K$-span may be identified with the commutative split semisimple $K$-algebra $K^{\weird C_R}$ (addition and multiplication in this algebra work component-wise). Similarly we may think of the set
$\Lambda_{LR} := e_L \Lambda e_R$ for $R,L\in\weird  R$ as sitting inside $K^{\weird C_R\cap \weird C_L}$. The set $\Lambda_{LR}$ 
may be construed as a $\Lambda_{LL}$-$\Lambda_{RR}$-bimdoule.
In short, in \cite{NebeSl2Char2} Nebe succeeds in describing the $\OO$-orders $\Lambda_{RR}$ and the sets
$\Lambda_{LR}$ as $\Lambda_{LL}$-$\Lambda_{RR}$-bimodules. However, the bimodule structure of 
$\Lambda_{LR}$ is not sufficient to describe $\Lambda$, since the multiplication maps $\Lambda_{LR}\times \Lambda_{RS}\longrightarrow \Lambda_{LS}$ cannot be fully recovered from the bimodule structure on the involved sets $\Lambda_{LR}$, $\Lambda_{RS}$ and $\Lambda_{LS}$. 

The first step in \cite{NebeSl2Char2} is to lift a $\bar k$-basis of $\bar{\Lambda}_{RR}$ to an $\OO$-basis of $\Lambda_{RR}$ (for each $R\in \weird R$). 
The $k$-basis used for this purpose was given in \cite{KoshitaChar2} as follows:
\begin{thm}[Koshita]
	Let $I\subset N$ and let $i \in N - I$. Let $j = j(i,I)$ be the unique integer $\leq i$ such that $j-1\notin I$ but $l\in I$ for all $j \leq l < i$. Define 
	\begin{equation}
		\omega_{i,I} := \alpha_{j, I} \cdot \alpha_{j+1} \cdots \alpha_{i-1} \cdot \alpha_{i}\cdot \alpha_i \cdot \alpha_{i-1} \cdots \alpha_{j+1} \cdot \alpha_j \in \bar \Lambda_{II}
	\end{equation}
	For a subset $T \subset N - I$ define 
	\begin{equation}
		\omega_{I,T}  := \prod_{i\in T} \omega_{I,i}  \in \bar \Lambda_{II}
	\end{equation}
	This product is well-defined independent of the order of the factors since $\bar \Lambda_{II}$ is commutative.
	The elements $\omega_{I,T}$ form a $\bar k$-basis of $\bar{\Lambda}_{II}$.
\end{thm}
Let $\widehat \alpha_{i,I} \in \Lambda_{I,I+\{i\}}$ be lifts of the elements $\alpha_{i,I}$.
One key observation in \cite{NebeSl2Char2} is that since each $\Lambda_{IJ}$ sits inside 
$K^{\weird C_I\cap \weird C_J}$ (which we may in turn view as a subset of $K^{\weird C}$ by simply extending vectors by zero)  
we can reorder elements in a product arbitrarily and always obtain the same result (this is only partially reflected in the commutativity relations in Koshita's presentation of $\bar \Lambda$, since we may also reorder the elements in a product in such a way that the start and endpoint of the corresponding path changes).
The reason is of course that the ring $K^{\mathcal C}$ (with component-wise multiplication) is commutative, and we may consider all products as being taken within this ring (we will do this frequently below).
So for instance $\widehat \alpha_{i,I}\cdot \widehat \alpha_{i, I+\{1\}}$ is equal to $\widehat \alpha_{i,I+\{i\}}\cdot \widehat \alpha_{i,I}$ inside $K^{\weird C}$. Now \cite[Lemma 3.10]{NebeSl2Char2} states that
$\frac{1}{2}\cdot \widehat \alpha_{i,I+\{i\}}\cdot \widehat \alpha_{i,I}$ lies in
$\Lambda_{I+\{i\}, I+\{i\}}$ (since $\alpha_{i,I+\{i\}}\cdot \alpha_{i,I}=0$ in $\bar{\Lambda}$), and is in fact a unit in this ring.
Let $u_{i,I}\in \Lambda_{I+\{i\}, I+\{i\}}$ denote its inverse. Then 
$u_{i,I} \cdot  \widehat \alpha_{i,I+\{i\}}\cdot \widehat \alpha_{i,I} = 2\cdot \eps_{{I+\{1\}}}$, where 
$\eps_{I+\{1\}}$ denotes the element in $K^{\weird C}$ which has entry equal to one in the components indexed 
by elements of $\weird C_{I+\{1\}}$, and entries equal to zero elsewhere. Since we may reorder elements in the product  we
obtain that $\widehat \alpha_{i,I} \cdot u_{i,I} \cdot  \widehat \alpha_{i,I+\{i\}} = 2\cdot \eps_{{I+\{1\}}}$ (note that this is now an element of $\Lambda_{I,I}$).
The same principle is applied to the elements $\omega_{I,i}$ defined above. First observe that
\begin{equation}
\widehat \alpha_{j, I} \cdot \widehat \alpha_{j+1} \cdots \widehat \alpha_{i-1} \cdot \widehat \alpha_{i} \cdot \widehat \alpha_{i} \cdot \widehat \alpha_{i-1} \cdots \widehat \alpha_{j+1} \cdot \widehat \alpha_j 
= (\widehat \alpha_{j, I}  \widehat \alpha_{j, I + \{j\}}) \cdots (\widehat \alpha_{i, I+\{j,\ldots,i-1\}}  \widehat \alpha_{i, I + \{j,\ldots,i\}})
\end{equation}
where the product on the right hand side is formed within $K^{\weird C}$. As we saw above, for each $j\leq l \leq i$ there is a unit $u_l$ in $\Lambda_{I+\{j,\ldots,l\}, I+\{j,\ldots,l\}}$ such that 
\begin{equation}
\widehat \alpha_{l, I+\{j,\ldots,l-1\}} \cdot u_l \cdot \widehat \alpha_{l, I+\{j,\ldots,l\}} = 2 \cdot \eps_{{I+\{j, \ldots, l\}}}
\end{equation}
We have hence found an explicit description of some element in $\Lambda_{I,I}$ which is analogous to the element
$\omega_{i,I}\in \bar\Lambda_{I,I}$ (however, it does not necessarily reduce to this element upon reduction modulo two):
\begin{equation}
	\beta_{i,I} := \widehat \alpha_{j, I}\cdot u_j \cdot \widehat \alpha_{j+1}\cdot u_{j+1} \cdots \widehat \alpha_{i-1}\cdot u_{i-1} \cdot \widehat \alpha_{i}\cdot u_i \cdot \widehat \alpha_{i-1}\cdot \widehat \alpha_{i-2} \cdots \widehat \alpha_{j+1} \cdot \widehat \alpha_j 
\end{equation}
By reordering the factors and using the definition of the $u_l$ one easily sees that
\begin{equation}\label{eqn_cgss2bnk3b}
	\beta_{i,I} = 2^{i-j+1}\cdot \eps_{I} \cdot  \eps_{I+\{j,\ldots,i\}} 
\end{equation}
\begin{thm}[{\cite[Theorem 3.12]{NebeSl2Char2}}]
	For any subset $I\subseteq N$ and any subset $T \subseteq N - I$ define
	\begin{equation}
		\beta_{T,I} := \prod_{i\in T} \beta_{I,i}
	\end{equation}
	where the empty product is defined to be $\eps_{I}$. 
	Then the $\beta_{T,I}$ form an $\OO$-basis of the $\OO$-order $\Lambda_{I,I}$. 
\end{thm}
Thanks to formula (\ref{eqn_cgss2bnk3b}) this description of $\Lambda_{I,I}$ is perfectly explicit.
Now let $I,J\subseteq N$ be two distinct subsets. Then we get the following information on
the $\Lambda_{I,J}$:
\begin{thm}[{\cite[Theorem 3.12]{NebeSl2Char2}}]
	If $\bar \Lambda_{I,J}\neq 0$ then
	\begin{equation}\label{eqn_sjhj6dghd}
		\Lambda_{I,J}  \iso \eps_I \cdot \Lambda_{I\cap J, I \cap J} \cdot  \eps_J
	\end{equation}
	as a $\Lambda_{I,I}$-$\Lambda_{J,J}$-bimdoule.
\end{thm}
For a full description of the order $\Lambda$, we need more than a bimodule-isomorphism in (\ref{eqn_sjhj6dghd}). In fact, (\ref{eqn_sjhj6dghd}) fixes 
$\Lambda_{I,J}$ exactly up to a $K\otimes\Lambda_{I,I}$-$K\otimes\Lambda_{J,J}$-bimodule-automorphism of $K\otimes \Lambda_{I,J}\iso K^{\weird C_I\cap \weird C_J}$. These bimodule automorphisms 
of $K^{\weird C_I\cap \weird C_J}$ may be identified with elements of $(K-\{0\})^{\weird C_I\cap \weird C_J}$ acting on $K^{\weird C_I\cap \weird C_J}$ by component-wise multiplication.
Thus, $\Lambda_{I,J}  \iso  \mu_{I,J}\cdot \eps_I\cdot \Lambda_{I\cap J, I \cap J} \cdot  \eps_J$ with $\mu_{I,J} \in (K-\{0\})^{\weird C_I\cap \weird C_J}$.
In \cite{NebeSl2Char2} the following information on $\mu_{I,J}$ is obtained (one should keep in mind though that the $\mu_{I,J}$ are not uniquely determined; the main source of the ambiguity is that the order $\Lambda$ is only well-defined up to conjugation)
\begin{thm}
	We may choose $\mu_{I,J}$ such that 
	\begin{equation}
		\mu_{I,J} = u_{I,J} \cdot 2^{|I-J|}\cdot \eps_{I}\cdot \eps_J
	\end{equation}
	where $u_{I,J}\in (\OO^\times)^{\weird C_I\cap \weird C_J}$.
\end{thm}
Conjecturally, we may choose all of the $u_{I,J}$ to be identical one. This would describe the order $\Lambda$ up to isomorphism. This order has (by construction)
semisimple $K$-span, the same decomposition matrix as the basic order of $\OO \SL_2(2^f)$ and it is self-dual with respect to the appropriate 
trace bilinear form. Moreover, it reduces to an $k$-algebra which, upon tensoring with $\bar k$, becomes isomorphic to the basic algebra  of $\bar k \SL_2(2^f)$ as described by Koshita.
In the present article we confirm this conjecture, as well as the analogous conjecture concerning the group ring of $\SL_2(p^f)$
proposed in \cite{NebeSl2Charp}.

\section{Transfer of unique lifting via derived equivalences}

In this section we cite the necessary theorems from \cite{EiseleDerEq}.
They establish the main technical tool used in this paper: a bijection between the sets of lifts (in the sense of the definition below) of two derived equivalent $k$-algebras. This bijection will allow us to shift the problem of proving that a given $k$-algebra lifts uniquely to an $\OO$-order to an analogous problem over a simpler algebra which is derived equivalent to the original one. 

\begin{defi}\label{defi_lifts}
	For a finite-dimensional $k$-algebra $\overline\Lambda$ define its \emph{set of lifts} as follows:
	\begin{equation}
		\Lifts(\overline{\Lambda}) := \left\{ (\Lambda, \varphi) \ | \  \textrm{$\Lambda$ is an $\OO$-order and $\varphi: k\otimes\Lambda \stackrel{\sim}{\rightarrow} \overline \Lambda$ is an isomorphism}\right\}
		\bigg / \sim
	\end{equation}
	where we say $(\Lambda, \varphi) \sim (\Lambda',\varphi')$ if and only if 
	\begin{enumerate}
	 \item There is an isomorphism $\alpha: \Lambda \stackrel{\sim}{\rightarrow} \Lambda'$
	 \item There is a $\beta \in \Aut_k(\overline\Lambda)$ such that the functor  $-\otimes^{\mathbb L}_{\overline\Lambda} {_\beta \overline\Lambda_\id}$ fixes all isomorphism classes of tilting complexes in $\mathcal K^b(\projC_{\overline\Lambda})$
	\end{enumerate}
	such that $\varphi=\beta\circ\varphi'\circ(\id_k\otimes \alpha)$.
	
	Moreover we define 
	\begin{equation}
	\mathfrak L(\bar \Lambda) := \{ \textrm{ Isomorphism classes of $\OO$-orders $\Lambda$ with $k\otimes \Lambda\iso \bar \Lambda$ }  \}
	\end{equation}
	and the projection map
	\begin{equation}
		\Pi:\ \Lifts(\bar \Lambda) \longrightarrow \mathfrak L(\bar{\Lambda})
	\end{equation}	
	
	Finally, we define the set of lifts with semisimple $K$-span
	\begin{equation}
		\Lifts_s(\overline\Lambda) := \{ (\Lambda,\varphi)\in\Lifts(\overline\Lambda) \ | \ K\otimes\Lambda \textrm{ is semisimple } \}	 
	\end{equation}
	and similarly	
	\begin{equation}
			\mathfrak L_s(\bar\Lambda) := \{ \Lambda\in \mathfrak L(\bar{\Lambda})\ | \ K\otimes\Lambda \textrm{ is semisimple } \}	 
	\end{equation}
\end{defi}

\begin{thm}[{\cite[Theorem 5.2]{EiseleDerEq}}]\label{thm_bijection_with_conditions}
	Let $\overline{\Lambda}$ and $\overline{\Gamma}$ be finite-dimensional $k$-algebras that are derived equivalent.
	Let the derived equivalence be afforded by the two-sided tilting complex 
	$X$. 
	Then there is a bijective map
	\begin{equation}
		\Phi_X:\ \Lifts(\bar \Lambda) \longrightarrow \Lifts(\bar \Gamma)
	\end{equation} 
	such that all of the following properties hold:
	\begin{enumerate} 
	\item[(i)] If $(\Lambda, \varphi)\in \Lifts(\bar{\Lambda})$ and $(\Gamma, \psi) = \Phi_X(\Lambda, \varphi)$, then
		there is a derived equivalence between $\Lambda$ and $\Gamma$.
	\item[(ii)] 
	 	$\Phi_X$ induces a bijection
		\begin{equation}
			\Lifts_s(\overline\Lambda) \longleftrightarrow \Lifts_s(\overline\Gamma)
		\end{equation}
	\item[(iii)] Set $\Phi := \Pi \circ \Phi_X$. If $(\Lambda,\varphi), (\Lambda',\varphi')\in \Lifts(\overline{\Lambda})$ are two lifts with $Z(K\otimes \Lambda) \iso Z(K\otimes \Lambda')$, then 
	\begin{equation}
		Z(K\otimes \Phi(\Lambda,\varphi)) \iso Z(K\otimes \Phi(\Lambda',\varphi'))
	\end{equation}
	and every choice of an isomorphism $\gamma: Z(K\otimes \Lambda)\rightarrow Z(K\otimes \Lambda')$ gives rise to a (canonically defined) isomorphism $\Phi(\gamma): Z(K\otimes \Phi(\Lambda,\varphi))\rightarrow Z(K\otimes \Phi(\Lambda',\varphi'))$.
	\item[(iv)] If $(\Lambda,\varphi), (\Lambda',\varphi')\in \Lifts(\overline{\Lambda})$ are two lifts and  $\gamma: Z( \Lambda) \stackrel{\sim}{\rightarrow} Z(\Lambda')$ is an isomorphism, then $\Phi(\gamma): Z(\Phi(\Lambda,\varphi)) \rightarrow Z(\Phi(\Lambda',\varphi'))$ is well-defined and an isomorphism as well.
	\item[(v)] If $(\Lambda,\varphi), (\Lambda',\varphi')\in \Lifts_s(\overline{\Lambda})$ are two lifts, and
	$\gamma: Z(K\otimes \Lambda) \stackrel{\sim}{\rightarrow} Z(K\otimes \Lambda')$ is an isomorphism such that  $D^\Lambda = D^{\Lambda'}$
	up to permutation of columns
	(where rows are identified via $\gamma$), then 
	$D^{\Phi(\Lambda,\varphi)}=D^{\Phi(\Lambda',\varphi')}$ up to permutation of columns (where rows are identified via $\Phi(\gamma)$).
	\item[(vi)]  If $(\Lambda,\varphi), (\Lambda',\varphi')\in \Lifts_s(\overline{\Lambda})$ are two lifts with  $D^\Lambda = D^{\Lambda'}$
	up to permutation of rows and columns then
	$D^{\Phi(\Lambda,\varphi)} = D^{\Phi(\Lambda',\varphi')}$
	up to permutation of rows and columns.
	\end{enumerate}
\end{thm}

\begin{thm}[{see \cite[Theorem 4.7]{EiseleDerEq}}]\label{thm_symm_der_eq}
	Let $\Lambda$ and $\Gamma$ be two derived-equivalent $\OO$-orders with semisimple $K$-span.
	Then we may identify $Z(K\otimes \Lambda)$ and $Z(K\otimes \Gamma)$. The order $\Lambda$ is self-dual 
	with respect to $T_u$ (with $u\in Z(K\otimes \Lambda)$) if and only if $\Gamma$ is self-dual with respect to 
	$T_{\tilde u}$, where $\tilde u \in Z(K\otimes \Gamma)$ is obtained from $u$ by flipping the signs 
	in some Wedderburn components (in \cite[Theorem 4.7]{EiseleDerEq} there is an explanation exactly which signs need flipping, 
	but this will not matter in the present paper).
	
	In the setting of the preceding theorem the following holds:
	Let $(\Lambda, \varphi)\in \Lifts(\bar{\Lambda})$ and $(\Gamma,\psi) := \Phi_X(\Lambda, \varphi)$.
	By the first point of the preceding theorem there is an isomorphism $\gamma: Z(K\otimes \Lambda) \rightarrow Z(K\otimes\Gamma)$.
 	Then 
	$\Lambda$ is self-dual with respect to $u\in Z(K\otimes\Lambda)$ if and only if 
	$\Gamma$ is self-dual with respect to $\tilde u \in Z(K\otimes \Gamma)$, where $\tilde u$ is obtained from $\gamma(u)$ by flipping signs in certain Wedderburn components. 
\end{thm}

We are actually interested in isomorphism classes of orders which reduce to a given $k$-algebra $\bar{\Lambda}$ 
, i.e. the set $\mathfrak L(\bar \Lambda)$. However, Theorem \ref{thm_bijection_with_conditions} only relates 
the sets $\Lifts(\bar{\Lambda})$ among derived equivalent algebras. Proposition \ref{prop_unique_lift_class_ext} below relates $\mathfrak L(\bar{\Lambda})$ and $\Lifts(\bar\Lambda)$ with each other in a special case (which will be sufficient for us).  It is a slightly strengthened version of \cite[Proposition 3.12]{EiseleDerEq} (strengthened in that it no longer requires that $k$ be algebraically closed).

\begin{prop}[{see \cite[Corollary 2.14]{EiseleDerEq}}]\label{prop_out0_triv_tilting}
	Assume $k$ is algebraically closed and let $A$ be a finite-dimensional $k$-algebra. Let $T\in \mathcal K^b(\projC_A)$ be a one-sided tilting complex. Then
	\begin{equation}
		T \otimes_A {_\id A_\gamma} \iso T \quad \textrm{ for all $\gamma \in \Out_k^0(A)$}
	\end{equation}
\end{prop}

\begin{prop}\label{prop_tilting_iso_extension_of_scalars}
    Let $A$ be a finite-dimensional $k$-algebra and let $S$ and $T$ be two tilting complexes over $A$. Then
    $S\iso T$ (in $\mathcal D^b(A)$) if and only if $\bar k\otimes S \iso \bar k \otimes T$ in
    $\mathcal D^b(\bar k\otimes A)$.
\begin{proof}
    If $\bar k \otimes S \iso \bar k\otimes T$, then there has to be some finite extension 
    $k'$ of $k$ such that $k'\otimes S \iso k'\otimes T$. By restriction we will also have
    $S^{[k':k]} \iso T^{[k':k]}$. There is a $k$-algebra $B$ and an invertible
    complex $X$ of $A$-$B$-bimodules such that
    $S\otimes_{A}^{\mathbb L} X$ is the stalk complex of a module. 
    But then $T\otimes_{A}^{\mathbb L} X$ will be the stalk complex of a module
    as well, since it becomes isomorphic to $S\otimes_{A}^{\mathbb L} X$ upon tensoring with $k'$
    (note that the functors $-\otimes^{\mathbb L}_{A} X$ and $k'\otimes_k-$ commute with each other; 
    also, tilting complexes which are stalk complexes of modules are distinguished by the fact that they have non-trivial homology in only a single degree).
    Now we can simply apply Krull-Schmidt theorem. So
    $S^{[k':k]}\otimes_{A}^{\mathbb L} X \iso T^{[k':k]}\otimes_{A}^{\mathbb L} X$
    implies that $S\otimes_{A}^{\mathbb L} X\iso T\otimes_{A}^{\mathbb L} X$ and therefore $S\iso T$.
\end{proof}
\end{prop}

Note that for any $k$-algebra $\bar \Lambda$ there is a left action of $\Out_k(\bar\Lambda)$ on $\Lifts(\bar{\Lambda})$. 
If $(\Lambda,\varphi)\in\Lifts(\bar{\Lambda})$ and $\alpha\in \Out_k(\bar{\Lambda})$ we simply set $\alpha \cdot (\Lambda, \varphi) := (\Lambda, \alpha\circ \varphi)$. It is proved in \cite[Proposition 3.7]{EiseleDerEq} that this is indeed well-defined (i. e. independent of the choice of a representative for $\alpha$).

\begin{corollary}\label{corollary_G_trivial}
	Let $\bar \Lambda$ be an finite-dimensional $k$-algebra, and let $G \leq \Out_k(\bar \Lambda)$ be a subgroup such that the $\bar k$-linear extensions of the elements of $G$ all lie in 
	$\Out_{\bar k}^0(\bar k \otimes \bar{\Lambda})$. Then $G$ acts trivially on $\Lifts (\bar \Lambda)$.
\begin{proof}
	Since $G$ acts trivially on isomorphism classes of tilting complexes in $\mathcal K^b(\projC_{\bar k \otimes \bar{\Lambda}})$ by Proposition \ref{prop_out0_triv_tilting}, it follows using 
	Proposition \ref{prop_tilting_iso_extension_of_scalars} that $G$ acts trivially on isomorphism classes of tilting complexes in
	$\mathcal K^b(\projC_{\bar{\Lambda}})$. But by definition of the equivalence relation ``$\sim$'' this means that $G$ acts trivially on $\Lifts(\bar{\Lambda})$. 
\end{proof}
\end{corollary}

\begin{prop}[{cf. \cite[Proposition 3.12]{EiseleDerEq}}]\label{prop_unique_lift_class_ext}
  Let $\Lambda \in \mathfrak L(\overline\Lambda)$, and let $\gamma: k\otimes\Lambda\stackrel{\sim}{\rightarrow}\overline\Lambda$.
    be an isomorphism. Now assume
   \begin{equation}\label{eqn_JHKHIOdjkswjdiew}
   	\overline{\Aut_\OO(\Lambda)} \cdot G = \Out_k(\overline\Lambda)
   \end{equation}
    where $\overline{\Aut_\OO(\Lambda)}$ is  the image of $\Aut_\OO(\Lambda)$ in $\Out_k(\overline\Lambda)$
        (here we identify $k\otimes\Lambda$ with $\overline\Lambda$ via $\gamma$)
    and $G\leq \Out_k(\overline\Lambda)$ is a subgroup such that the $\bar k$-linear  extensions of all elements of $G$ lie in $\Out_{\bar k}^0(\bar k \otimes_k \overline \Lambda)$.
    Then the fiber $\Pi^{-1}(\{\Lambda\})$ has cardinality one.    
\begin{proof}
        Let  $(\Lambda, \varphi)\in\Lifts(\overline\Lambda)$ for some $\varphi: k\otimes\Lambda\stackrel{\sim}{\longrightarrow}\bar\Lambda$ (i. e. $(\Lambda,\varphi)$ is an arbitrary element in $\Pi^{-1}(\{\Lambda\})$).
        We intend to show $(\Lambda, \varphi)\sim (\Lambda, \gamma)$, since this will imply that $\Pi^{-1}(\{\Lambda\})$ contains indeed only a single element.
        Now if (\ref{eqn_JHKHIOdjkswjdiew}) holds, we can write $\gamma\circ\varphi^{-1}=\gamma\circ (\id_k\otimes \hat\alpha) \circ \gamma^{-1}\circ \beta$ for some $\hat\alpha\in \Aut_\OO(\Lambda)$ and $\beta \in G$. Hence $\gamma \circ (\id_k\otimes\hat\alpha^{-1})=  \beta \circ \varphi$.
        Corollary \ref{corollary_G_trivial} (together with the definition of ``$\sim$'') implies $(\Lambda,\gamma)\sim (\Lambda,\beta^{-1}\circ \gamma\circ(\id_k\otimes\hat\alpha^{-1}))=(\Lambda, \varphi)$.
\end{proof}
\end{prop}

\section{The algebra $k\Delta_2(p^f)$ and unique lifting}

In this section we will write $k\Delta_2(p^f)$ explicitly as a quotient of a quiver algebra
(at least in the case when $k$ splits $\Delta_2(p^f)$),
and use this presentation to show that it lifts uniquely to an $\OO$-order satisfying certain properties. At least the first part of this (finding a presentation as a quotient of a quiver algebra) is relatively straightforward.
The reason for looking at the group algebra of $\Delta_2(p^f)$ is that its blocks (one block if $p=2$, two blocks otherwise)
are the Brauer correspondents of the blocks of maximal defect of the group algebra of $\SL_2(p^f)$. Other than those blocks of maximal defect, the group algebra 
of $\SL_2(p^f)$ only has a block of defect zero. This block of defect zero will not be of interest to us though, since all questions we are concerned with can be answered trivially for such a block (after all, a block of defect zero is just a matrix ring).

\begin{defi}
	Assume that $A$ is an abelian $p'$-group such that $kA$ is split. Denote by $\hat{A}$ the character group of $A$, that is, $\Hom(A,k^\times)$
	(abstractly we will have $A \iso \hat{A}$). Assume moreover that $A$ is acting on a $p$-group $P$ by automorphisms. Let 
	\begin{equation}
		\Jac(kP)/\Jac^2(kP) \iso \bigoplus_{i=1}^l S_i
	\end{equation}
	be a decomposition of $\Jac(kP)/\Jac^2(kP)$ as a direct sum of simple $kA$-modules $S_1,\ldots,S_l$.
	We define the set $X(P,A)$ to be the disjoint union
	\begin{equation}
		\biguplus_{i=1}^l \{ \chi_{S_i} \}
	\end{equation}
	where $\chi_{S_i} \in \hat A$ denotes the character of $A$ associated to $S_i$.
\end{defi}

\begin{lemma}\label{lemma_action_radical}
	Let $P = C_p^f$ and let $A$ be a group acting on $P$ by automorphisms. View $P$ as
	an $\F_p$-vector space by identifying $C_p^f$ with $(\F_p^f,+)$.  
	Under this identification, $P$ becomes an $\F_pA$ module. 
	Then
	\begin{equation}
		\Jac(kP)/\Jac^2(kP) \iso_{kA} k\otimes_{\F_p} P
	\end{equation}
\begin{proof}
	First note that after identifying $P$ with $\F_p^f$, the fact that $A$ acts on $P$ by automorphisms
	translates into $A$ acting linearly on $\F_p^f$, as each automorphism of $(\F_p^f,+)$ is automatically
	$\F_p$-linear. This turns $P$ into an $\F_pA$-module (in fact, the isomorphism type of this module 
	is independent of the choice of the  identification of $P$ with $\F_p^f$).
	Let $x_1,\ldots,x_f$ be a minimal generating system for $P=C_p^f$. Then
	$1\otimes x_1,\ldots,1\otimes x_f$ is a $k$-basis for $k\otimes_{\F_p} P$.
	Now define a $k$-linear map
	\begin{equation}
		\Phi:\ k\otimes_{\F_p} P \rightarrow \Jac(kP)/\Jac^2(kP):\ 1\otimes x_i \mapsto x_i - 1
	\end{equation}
	Since the $x_i-1$ lie in $\Jac(kP)$ and they are a minimal (with respect to inclusion) generating set for $kP$ as a $k$-algebra, they
	form a $k=kP/\Jac(kP)$ basis of $\Jac(kP)/\Jac^2(kP)$. Hence $\Phi$ is an isomorphism of vector spaces.
	We only need to check that $\Phi$ is $A$-equivariant (or, more generally, $\Aut(P)$-equivariant). This amounts to showing that for all
	$n_1,\ldots,n_f \in \Z_{\geq 0}$ the following holds:
	\begin{equation}\label{eqn_dhjiowhjiohjerioe}
		x_1^{n_1}\cdots x_f^{n_f} - 1 \equiv \sum_{i=1}^f n_i\cdot (x_i-1) \mod \Jac^2(kP) 		
	\end{equation}
	Let $x,y \in P$. Then clearly $(x-1)(y-1) \in \Jac^2(P)$, and hence $xy-x-y+1 \equiv 0 \mod \Jac^2(kP)$. This can be rewritten as $xy-1 \equiv (x-1)+(y-1) \mod \Jac^2(kP)$. Applying this equality iteratedly clearly implies (\ref{eqn_dhjiowhjiohjerioe}).
\end{proof}
\end{lemma}

\begin{prop}\label{prop_rel_delta22f}
	Let $G=P \rtimes A$ with $P \iso C_p^f$ and $A$ an abelian $p'$-group acting on $P$. 
	If $k$ splits $G$ then 
	\begin{equation}
		kG \iso kQ/I
	\end{equation}
	where $Q$ is the quiver which has vertices $e_{\chi}$ in bijection with the elements $\chi \in \hat{A}$, 
	and an arrow $e_\chi \stackrel{s_{\chi,\psi}}{\longrightarrow} e_{\chi\cdot \psi}$ for each $\chi \in \hat{A}$
	and $\psi\in X(P,A)$. $I$ is the ideal generated by the relations 
	\begin{equation}\label{rel_ext_one}
		s_{\chi, \psi} \cdot s_{\chi\cdot \psi, \varphi}=s_{\chi, \varphi} \cdot s_{\chi\cdot \varphi, \psi}	 
		\quad\textrm{ for all $\chi \in \hat A$ and $\psi,\varphi \in X(P,A)$}
	\end{equation}
	and
	\begin{equation}\label{rel_ext_two}
		\prod_{i=0}^{p-1} s_{\chi\cdot \psi^i,\psi} = 0 \quad \textrm{ for all $\chi \in \hat A$ and $\psi \in X(P,A)$}
	\end{equation}
\begin{proof}
	We first look at $kP$. We have $kC_p \iso k[T]/\langle T^p \rangle$, and 
	\begin{equation}
		kP \iso \bigotimes^f kC_p \iso k[T_1,\ldots,T_{f}]/( T_1^p,\ldots, T_{f}^p )
	\end{equation}
	Given any minimal generating set 
	$t_1,\ldots,t_{f}$ of $kP$
	contained in $\Jac(kP)$, the epimorphism $k[T_1,\ldots,T_{f}]\twoheadrightarrow kP$
	sending $T_i$ to $t_i$ has the same kernel $(T_1^p,\ldots,T_{f}^p )$. 
	This is simply because any automorphism of $k[T_1,\ldots,T_f]$ mapping the ideal 
	$(T_1,\ldots,T_f)$ into itself will map the ideal $(T_1^p,\ldots,T_f^p)$ into itself as well.

	Now consider the action of $A$ on $\Jac(kP)$ by conjugation. Since $kA$ is abelian and split semisimple, there is
	a basis $t_1,\ldots,t_{p^f-1}$ of $\Jac(kP)$ such that for each $i$ the conjugates $u^{-1}t_iu$ are a multiple of $t_i$ for all $u\in A$. We may choose a minimal generating set for $kP$ from said $t_i$'s, say 
	(after reindexing) $t_1,\ldots,t_f$. As the images of $t_1,\ldots,t_f$ in $\Jac(kP)/\Jac^2(kP)$ form a basis,
	there is a bijective map 
	\begin{equation}
	    X(P,A) \longrightarrow \{t_1,\ldots,t_f\}: \ \psi \mapsto s_{\psi}
	\end{equation}
	  such that
	$u^{-1}\cdot  s_\psi\cdot  u = \psi(u)\cdot s_{\psi}$ for all $u \in A$.
	Define furthermore for each $\chi \in \hat{A}$ the corresponding primitive idempotent $e_{\chi} \in kA$
	via the standard formula
	\begin{equation}
		e_\chi = \frac{1}{|A|} \sum_{a \in A} \chi(a)\cdot a^{-1}
	\end{equation}
	This is a full set of orthogonal primitive idempotents in $kG$. Furthermore
	\begin{equation}
		e_\chi \cdot s_\psi = \frac{1}{|A|} \sum_{a \in A} \chi(a)\cdot a^{-1} s_\psi \cdot a \cdot a^{-1}
		= s_{\psi}\cdot \frac{1}{|A|} \sum_{a \in A} \chi(a)\psi(a)\cdot a^{-1}= s_{\psi} \cdot e_{\chi\cdot \psi}
	\end{equation}
	Hence define
	\begin{equation}
		s_{\chi,\psi} := e_{\chi} \cdot s_{\psi} \quad \textrm{ for all $\chi \in \hat{A}, \psi \in X(P,A)$}
	\end{equation}
	The fact that the $s_\psi$ commute implies the relation (\ref{rel_ext_one}), and the
	fact that $s_\psi^p=0$ implies relation (\ref{rel_ext_two}).
	What we have to verify though is that that the $s_{\psi}$ and $e_{\chi}$ generate $kG$ as a $k$-algebra, and that there
	are no further relations (i. e. $\dim_k kG = \dim_k kQ/I$).

% 	Now it is an easy argument (see for instance \cite[Lemma 5.8]{Alperin}) that
% 	$kG \cdot \Jac^i(kP) = \Jac^i(kG)$. This clearly implies 
% 	$kA \cdot \Jac^i(kP) = \Jac^i(kG)$ as well (as $kG = kA \cdot kP = kP\cdot kA$). 
	The $s_{\psi}$ generate $kP$ as a $k$-algebra and the $e_\chi$ generate $kA$ even as a $k$-vector space.
	Hence together they generate $kP \cdot kA = kG$ as a $k$-algebra. Now to the dimension of $kQ/I$. 
	We can use relation (\ref{rel_ext_one}) to rewrite a path involving the arrows 
	$s_{\chi_1,\psi_1},\ldots, s_{\chi_l,\psi_l}$ (in that order) as a path
	$s_{\tilde \chi_1,\tilde\psi_1}\cdots s_{\tilde \chi_l,\tilde\psi_l}$ for any chosen
	reordering
	$(\tilde\psi_1,\ldots,\tilde\psi_l)$ 
	of
	$(\psi_1,\ldots,\psi_l)$. Notice that necessarily $\chi_1=\tilde\chi_1$, and all other $\tilde\chi_i$ are
	determined by $\tilde\chi_1$ and the $\tilde\psi_i$. Also we may assume, due to relation (\ref{rel_ext_two}),
	that no $p$ of the $\psi_i$ are equal. So ultimately, there are at most $|\hat{A}|\cdot p^{|X(P,A)|}$
	linearly independent paths ($|\hat{A}|$ choices for the starting point $\chi_1$, $p$ choices for
	the number of occurrences of each element of $X(P,A)$ in the sequence $(\psi_1,\ldots,\psi_l)$).
	Hence
	\begin{equation}
		\dim kQ/I \leq |\hat{A}|\cdot p^{|X(P,A)|} = |A|\cdot p^f = \dim_k kG
	\end{equation}
	and thus the epimorphism $kQ/I \twoheadrightarrow kG$ is in fact an isomorphism.
\end{proof}
\end{prop}

\begin{remark}
	It seems practical to keep on using the notation
	\begin{equation}
		s_\psi = \sum_{\chi \in \hat{A}} s_{\chi,\psi}
	\end{equation}
	With this notation we may just write
	\begin{equation}
		kG \iso kQ/\left\langle  s_{\psi}s_{\varphi}-s_{\varphi}s_{\psi},\ s_{\psi}^p \ | \ \psi,\varphi \in X(P,A )\right \rangle
	\end{equation}
\end{remark}

\begin{prop}
	Let $G=\Delta_2(p^f)$, $P=\mathbb G_a(\F_{p^f})\iso C_p^f$ and $A=\mathbb{G}_m(\F_{p^f})\iso C_{p^f-1}$ (we view $P$ as the subgroup of $G$ consisting
	of diagonal matrices and $A$ as the subgroup of $G$ consisting of unipotent matrices).
	Assume $\F_{p^f} \subseteq k$ and identify $\hat{A} = \Z/(p^f-1)\Z$ (where we identify $i$ with the character that sends $a\in A$ to $a^i\in k^\times$) and write the group operation in $\hat{A}$ additively.
	 Then
	\begin{equation}
		X(P,A) = \left\{ 2\cdot p^q \ | \ q =0,\ldots,f-1 \right\}
	\end{equation}
	In particular, the $\Ext$-quiver $Q$ of $k\Delta_2(p^f)$ has $p^f-1$ vertices $e_i$ labeled by elements
	$i \in \Z/(p^f-1)\Z$. There are precisely $f$ arrows $s_{i,2\cdot p^q}$ (for $q\in\{0,\ldots,f-1\}$) emanating from each vertex $e_i$.
\begin{proof}
	$G=P\rtimes A$ is a semidirect product.
	The action of $A$ on $P$ is given by
	\begin{equation}
		P \times A \rightarrow P: \ (b,a) \mapsto b \cdot a^2 \quad \textrm{where we identified $A=\F_{p^f}^\times$, $P=\F_{p^f}$}
	\end{equation}
	Let us denote the $\F_pA$ module $\F_{p^f}$ with the action of $A$ specified above by $M$.
	According to Lemma \ref{lemma_action_radical} we have to determine the simple constituents
	of $k\otimes_{\F_p}M$ as a $k A$-module. Note that there is a (one-dimensional) $\F_{p^f}A$-module
	$\tilde M$ with $\tilde M|_{\F_pA} \iso M$. So clearly
	\begin{equation}
		k\otimes_{\F_p} M \iso \bigoplus_{\gamma\in\Gal(\F_{p^f} / \F_p)} k \otimes_{\F_{p^f}} \tilde M^\gamma 
	\end{equation}
	Now $\Gal(\F_{p^f} / \F_p) \iso C_f$ is generated by the Frobenius automorphism. So the simple
	constituents of $k\otimes_{\F_p} M$ are just copies of $k$ on which $a\in A$ acts as $a^{2\cdot p^q}$ for
	$q\in\{0,\ldots, f-1\}$. This shows that $X(P,A)$ is as claimed. The shape of the $\Ext$-quiver is now
	immediate from Lemma \ref{lemma_action_radical}.
\end{proof}
\end{prop}

\newcommand{\ul}[1]{[ \mathbf{#1} ]}

\begin{notation}
    We define symbols
    \begin{equation}
	   \ul q := 2\cdot p^q  
    \end{equation}
    to refer to the elements of $X(P,A)$ in the situation of the above proposition.
\end{notation}

\begin{lemma}
	Assume $k$ splits $\Delta_2(p^f)$.
	$k\Delta_2(p^f)$ consists of a single block if $p=2$, and two isomorphic blocks otherwise.
	In the case $p=2$, the Cartan matrix is given by $I+J$, where $I$ is the identity matrix, and
	$J$ is the matrix that has all entries equal to one. In the case $p$ odd, the Cartan matrix of either one of the
	two blocks is $I+2\cdot J$.
\begin{proof}
	The $(i,j)$-entry of the Cartan matrix is, by definition, the $k$-dimension of $e_i \cdot kQ/I\cdot e_j$.
	Let $E = \langle e_1,\ldots,e_{p^f-1}\rangle_k$ be the subspace of $kQ/I$ spanned by the idempotents.
	Clearly, $kQ/I = E\oplus \Rad(kQ/I)$. So $\dim_k e_i \cdot kQ/I\cdot e_j = \delta_{ij} + \dim_k e_i \Rad(kQ/I) e_j$. Now,
	using the quiver relations from Proposition \ref{prop_rel_delta22f}, we can deduce that $\dim_k e_i \Rad(kQ/I) e_j$
	is equal to the number of vectors $(0,\ldots,0)\neq (n_0,\ldots,n_{f-1})\in \{0,\ldots,p-1\}^f$ such that
	\begin{equation}\label{eqn_iouhjeoi}
		2\cdot \sum_{q=0}^{f-1} n_q \cdot p^q \equiv i-j \mod (p^f-1)
	\end{equation}
	If $p$ is odd and $i-j$ is odd as well, then (since $p^f-1$ will be even) the congruence cannot possibly be satisfied by any sequence of $n_q$'s. So the corresponding entries in the Cartan matrix are zero. 
	Now assume that $p$ is odd and $i-j$ is even. Then the above congruence is equivalent to
	\begin{equation}
		\sum_{q=0}^{f-1} n_q \cdot p^q \equiv \frac{i-j}{2} \mod \left(\frac{p^f-1}{2} \right)
	\end{equation}
	By uniqueness of the $p$-adic expansion of an integer, the analogous equation modulo $p^f-1$ has a 
	unique solution (in the case $i-j\equiv 0 \mod (p^f-1)$ we would have two solutions, but we said above that we only consider solutions where not all of the  $n_q$'s are zero). Hence the equation above has precisely two solutions. 

	Now if $p=2$, the factor ``$2$'' in (\ref{eqn_iouhjeoi}) is a unit in the ring $\Z/(2^f-1)\Z$, and hence can be divided out. The remaining equation has a unique solution thanks to the uniqueness of the $2$-adic expansion of
	an integer (again discounting the zero solution). 
\end{proof}
\end{lemma}

\begin{remark}
	By counting conjugacy classes in the group $\Delta_2(2^f)$, one easily obtains that
	\begin{equation}
	  \dim_K Z(K\Delta_2(2^f))=2^f
	\end{equation}
	In the same way one obtains for $p$ odd that 
	\begin{equation}
	    \dim_K Z(K\Delta_2(p^f))=p^f+3
	\end{equation}
	Since $k\Delta_2(p^f)$ is the direct sum of two
	isomorphic blocks, the dimension of the center of either one of these blocks
	is $(p^f+3)/2$.
\end{remark}

For reasons that will become apparent in the section on descent to smaller fields, we would like to investigate 
a slightly larger class of algebras than the blocks of $k\Delta_2(p^f)$, namely those (split) $k$-algebras which 
become isomorphic to $k\Delta_2(p^f)$ upon extension of the 
ground field.

\begin{defi}
	We call a split $k$-algebra $\overline\Lambda$ with $\bar k\otimes \overline\Lambda \iso B_0(\bar k\Delta_2(p^f))$ a \emph{split $k$-form}
	of the principal block $B_0(\bar k\Delta_2(p^f))$ of $\bar k\Delta_2(p^f)$.
\end{defi}

\begin{remark}
	 If $\overline\Lambda$ is a split $k$-form of $B_0(\bar k \Delta_2(p^f))$, then $\overline\Lambda$ has the same
	 $\Ext$-quiver and the same Cartan matrix as $B_0(\bar k\Delta_2(p^f))$. 
	 Moreover, the $k$-dimension of the center of $\overline\Lambda$ is equal to
	 the $\bar k$-dimension of the center of $B_0(\bar k \Delta_2(p^f))$. 
\end{remark}

\begin{remark}
	 The quiver relations given in (\ref{rel_ext_one}) and (\ref{rel_ext_two}) are defined over $\F_p$. In particular, 
	 even if $k$ is no splitting field for $\Delta_2(p^f)$, the blocks of $kQ/I$ are split $k$-forms
	 of $B_0(\bar k \Delta_2(p^f))$.
\end{remark}

\begin{prop}[Shape of split $k$-forms]\label{prop_kform_delta2_rels}
    Let $\overline\Lambda$ be a split $k$-form of $B_0(\bar k\Delta_2(p^f))$. 
    By $Q$ we now denote the $\Ext$-quiver of $B_0(\bar k \Delta_2(p^f))$ (as opposed to the entire group ring $\bar k\Delta_2(p^f)$, which it was before).
    Denote (as before) the vertices of $Q$ by $e_{2i}$ and 
    the arrows by $s_{2i,q}$. Then $\overline \Lambda$ is isomorphic to $kQ/I'$ for some ideal $I'$ which contains relations
    \begin{equation}\label{eqn_jksjuu22nbh37}
	 \prod_{j=0}^{p-1} s_{2i+j\cdot\ul q,q} \quad \textrm{ for all $i\in \Z$ and $q\in \{0,\ldots,f-1\}$}
    \end{equation}
    and relations of the shape
    \begin{equation}\label{eqn_jk877bbsgz33hh}
	s_{2i,q} \cdot s_{2i+\ul q,q'} - \alpha_{2i,q,q'} \cdot s_{2i,q'}\cdot s_{2i+\ul {q'}, q}
    \end{equation}
    with $i$ ranging over $\Z$, $q$ and $q'$ ranging over $\{0,\ldots,f-1\}$ and the $\alpha_{2i,q,q'}$ being of the form 
    \begin{equation}c_{2i,q,q'}\cdot e_{2i} + r_{2i,q,q'}\end{equation}
    for some $c_{2i,q,q'}\in k^\times$ and some $k$-linear combination $r_{2i,q,q'}$ of closed paths of positive length starting and ending in $e_{2i}$ 
    (hence, by construction, the $\alpha_{2i,q,q'}$ will lie in $(e_{2i}\cdot kQ/I'\cdot e_{2i})^\times$). 
    
    The relations given in (\ref{eqn_jksjuu22nbh37}) and (\ref{eqn_jk877bbsgz33hh}) together with all
    paths of length $|\Delta_2(p^f)|$ (or any other sufficiently large number) generate $I'$.
\begin{proof}
    We can assume that $\overline\Lambda \iso kQ/I'$ for some ideal $I'$ contained in the ideal of $kQ$ generated by the paths of length at least two.
    We proceed to show that $I'$ is of the desired form.
    Choose an embedding $\varphi:\ kQ/I' \hookrightarrow \bar k Q/I$ that maps the idempotents $e_{2i}$ to themselves such that the $\bar k$-span
    of the image of $\varphi$ is all of $\bar k Q/I$. Then for each $i$ and $q$ the image $\varphi (s_{2i,q})$ has to be equal to $x_{2i,q}\cdot s_{2i,q}$ for
    some $x_{2i,q}\in (e_{2i}\cdot  \bar k Q/I \cdot e_{2i})^\times$ (since the relations in $I$ can be used to show that $e_{2i}\cdot \bar kQ/I\cdot e_{2i+\ul q}=e_{2i}\cdot \bar kQ/I\cdot e_{2i}\cdot s_{2i,q}$; now if $x_{2i,q}$ were no unit in $e_{2i}\cdot \bar kQ/I\cdot e_{2i}$, then $\varphi(s_{2i,q})$ would be contained in $\Jac^2(\bar kQ/I)$ and
    therefore the $\varphi(s_{2i,q})$ together with the $e_{2i}$ could not possibly generate $\bar kQ/I$ as a $\bar k$-algebra).
    Since the relations in $I$ imply that $e_{2i}\cdot \bar kQ/I \cdot e_{2i}\cdot s_{2i,q} = s_{2i,q} \cdot e_{2i+\ul q}\cdot \bar kQ/I \cdot e_{2i+\ul q}$, the relations in (\ref{eqn_jksjuu22nbh37}) follow immediately from the corresponding relation in $I$ by application of $\varphi$.
    
    Analogous to the above discussion, we can also deduce that for all $i\in\Z$ and $q,q'\in\{0,\ldots,f-1\}$
    \begin{equation}\label{eqn_kks883bbdggvgdzgzegzi3}
	\varphi(s_{2i,q})\cdot \varphi(s_{2i+\ul q, q'}) = \beta_{2i,q,q'} \cdot \varphi(s_{2i,q'})\cdot \varphi(s_{2i+\ul {q'}, q}) 
    \end{equation}
    for some $\beta_{2i,q,q'}\in (e_{2i}\cdot \bar kQ/I\cdot e_{2i})^\times$. Now take $\alpha'_{2i,q,q'} := (\id_{\bar k}\otimes_k \varphi)^{-1}(\beta_{2i,q,q'})\in \bar k\otimes_k k  Q/I'$. Choose a $k$-vector space complement $V$ of $k$ in $\bar k$ and choose $\alpha_{2i,q,q'}\in e_{2i} \cdot kQ/I'\cdot e_{2i}$ such that 
    $\alpha'_{2i,q,q'} =\alpha_{2i,q,q'}+\textrm{ (Sum of paths with coefficients in $V$) }$. Now clearly the following holds:
    \begin{equation}
	s_{2i,q}\cdot s_{2i+\ul q, q'} = \alpha_{2i,q,q'} \cdot s_{2i,q'}\cdot s_{2i+\ul {q'}, q} + \textrm{ (Sum of paths with coefficients in $V$)} 
    \end{equation}
      in $\bar k\otimes_k kQ/I'$. Since a sum of paths with coefficients in $V$ must be $k$-linearly independent from $kQ/I'$, the relation (\ref{eqn_jk877bbsgz33hh}) must
      hold with this choice of $\alpha_{2i,q,q'}$. To see that the coefficient of $e_{2i}$ in $\alpha_{2i,q,q'}$ is non-zero we could simply map the relation
      back into $\bar kQ/I$ using $\varphi$ and subtract it from relation (\ref{eqn_kks883bbdggvgdzgzegzi3}). This implies $(\beta_{2i,q,q'}-\varphi(\alpha_{2i,q,q'}))\cdot s_{2i,q'}\cdot s_{2i+\ul{q'},q}=0$, and hence $\beta_{2i,q,q'}-\varphi(\alpha_{2i,q,q'})$ is no unit in $e_{2i}\cdot \bar kQ/I\cdot e_{2i}$, which forces
      $\varphi(\alpha_{2i,q,q'})$ to be a unit.
      
      The claim that the given relations together with all paths of some sufficiently large length generate $I'$ can be verified by showing that they can be used to rewrite any path as a linear combination of paths of the form
      \begin{equation}\label{eqn_dhuhdkuhu3393u93dhijhjdkh}
      		s_{2i,q_1}\cdot s_{2i+\ul {q_1}, q_2}\cdots s_{2i+\ul{q_1}+\ldots+\ul{q_{l-1}},q_l} 
      \end{equation}
      such that $q_1\leq q_2 \leq \ldots \leq q_l$ and no $p$ of the $q_j$'s are equal. The latter requirement can be met using relation (\ref{eqn_jksjuu22nbh37}). If the $q_j$'s are not ordered as wanted, relation (\ref{eqn_jk877bbsgz33hh}) can be used to permute them. This will however produce some summands of strictly greater length. So one can apply a rewriting strategy where one starts with the paths of smallest length
      which are not already in the desired standard form, rewrites those (possibly altering or adding some summands of strictly  greater length) and then
      repeats the process until the shortest paths not in standard form are bigger than the cut-off length and therefore equal to zero.
\end{proof}
\end{prop}

\begin{lemma}\label{lemma_decomp_ksl22f}
Let $\overline\Lambda$ be a split $k$-form of $B_0(\bar k\Delta_2(p^f))$
\begin{enumerate}
\item Assume $p=2$. Then any lift $\Lambda\in \mathfrak{L}_s(\overline\Lambda)$ with
$\dim_K Z(K\otimes \Lambda)=\dim_k Z(\overline\Lambda)$ has the following decomposition matrix over a splitting field
\begin{equation}\label{decomp_kdelta2}
	\left[
		\begin{array}{cccc}
			1 & 0 & \cdots & 0 \\
			0 & 1 & \cdots & 0 \\
			\vdots & \vdots & \ddots & \vdots \\
			0 & 0 & \cdots & 1\\
			1 & 1 & \cdots & 1
		\end{array}
	\right]
\end{equation}
up to permutation of rows.
\item Assume $p\neq 2$. If $\Lambda\in \mathfrak{L}_s(\overline\Lambda)$ with
$\dim_K Z(K\otimes \Lambda)=\dim_k Z(\overline\Lambda)$ , then the
	decomposition matrix of $\Lambda$ over a splitting field  looks as follows:
	\begin{equation}\label{decomp_kdelta2p}
		\left[ \begin{array}{cccc} 1 & 0 & \cdots & 0 \\
			0 & 1 & \cdots & 0 \\
			\vdots &\vdots & \ddots & \vdots \\
			0 & 0 & \cdots & 1 \\
			1 & 1 & \cdots & 1 \\
			1 & 1 & \cdots & 1
		       \end{array}
		\right]
	\end{equation}
	up to permutation of rows.
\item 	Fix a $\Lambda\in \mathfrak{L}_s(\overline\Lambda)$ subject to the condition on the center
	as above. Assume that there is some totally ramified extension of $K$ that splits $\Lambda$.
	\begin{enumerate}
	\item If $p=2$, then $K$ already splits $\Lambda$.  
	\item If $p$ is odd then all one-dimensional representations of $\bar K\otimes \Lambda$ are already 
	defined over $K\otimes \Lambda$. If $K$ does not split $K\otimes \Lambda$, then $K\otimes \Lambda$ has
	a unique representation of dimension greater than one, and its endomorphism ring is a totally ramified
	extension of $K$ of degree two. In particular, in that case, the decomposition matrix of $\Lambda$ is as
	in (\ref{decomp_kdelta2p}) with the last row removed.
	\end{enumerate}
\end{enumerate}
\begin{proof}
	Concerning the first two parts:
	Let $D$ be the decomposition matrix of $\Lambda$ (over a splitting system). First note that all entries
	of $D$ must be $\leq 1$, as $D^\top\cdot D$ is equal to the Cartan matrix $C$ of $\bar k\Delta_2(p^f)$, which has ``$2$'''s (respectively ``$3$'''s) on the diagonal. It is straightforward to prove that the only solutions (with non-negative integer entries $\leq 1$) to the equation $D^\top \cdot D = C$ are, up to permutation of rows and columns, the ones given in statement of this lemma.

	Now we have a look at the assertions in the non-splitting case. 
	First assume that there is a simple $K\otimes \Lambda$-module $V$ such that $\End_{K\otimes\Lambda}(V)$ is non-commutative.
    Let $P$ be a projective indecomposable $\Lambda$-lattice (note that $k\otimes\Lambda\iso\overline\Lambda$ is split, so
    indecomposable projectives are absolutely indecomposable) such that $V$ occurs as a composition factor of $K\otimes P$.
    Since the endomorphism ring of $V$ is non-commutative, $\bar K\otimes V$ is not multiplicity-free, but it is still
    a composition factor of $\bar K\otimes P$. Hence there is some simple $\bar K \otimes \Lambda$-module
    which occurs in $\bar K\otimes P$ with multiplicity greater than one. This is the same as saying that (over a splitting system) there is a 
    decomposition number greater than one, which, as we have seen above, is impossible.
    Now let $V$ be any simple $K\otimes\Lambda$-module. As we have seen $E := \End_{K\otimes\Lambda}(V)$ is commutative, and therefore it is necessarily contained in any splitting field for
    $K\otimes\Lambda$. Since by assumption there is a splitting field that is totally ramified over $K$,
    the field $E$ must be totally ramified over $K$ as well.
	Now we look at how the decomposition matrix 
	over $K$ relates to the decomposition matrix over a splitting field.
	$\End_{\bar K\otimes\Lambda}(\bar K\otimes V) \iso \bar K \otimes_K E \iso \bigoplus^{\dim_K E} \bar{K}$. This implies that $\bar{K}\otimes V$ decomposes into $e := \dim_K E$
	non-isomorphic absolutely irreducible modules $V_1,\ldots, V_e$. Whenever $P$ is a projective indecomposable $\Lambda$-module, 
	the multiplicity of any $V_i$ in $\bar K\otimes P$ is the same as the multiplicity of $V$ in $K\otimes P$. Hence,
	the decomposition matrix of $\Lambda$ over a splitting field arises from the decomposition matrix 
	over $K$ by repeating certain rows. The shape of the decomposition matrix over a splitting field proved above then limits the simple $K\otimes \Lambda$-modules that may not be split sufficiently so that our claims follow.
\end{proof}
\end{lemma}

\begin{notation}\label{notation_cpi_ordering_sl2}
    Let $\Lambda$ be an $\OO$-order with semisimple $K$-span and let $\eps_1,\ldots,\eps_n\in Z(K\otimes\Lambda)$ be the
    central primitive idempotents. So, in particular, we have fixed a bijection $\{1,\ldots,n\}\leftrightarrow \{\textrm{ central primitive idempotents }\}$.
    \begin{enumerate}
    \item 
    Given an element $u\in Z(K\otimes\Lambda)$ we set 
    \begin{equation}
	u_i := \eps_i \cdot u \quad \textrm{for all $i\in\{1,\ldots,n\}$}
    \end{equation}
    \item When dealing with orders $\Lambda$ which have a decomposition matrix like the one in (\ref{decomp_kdelta2}) or (\ref{decomp_kdelta2p}), we make the following convention 
    concerning the ordering of the central primitive idempotents:
    We choose indices so that the idempotents associated to rows in the decomposition matrix with more than one
    non-zero entry come last.
    \end{enumerate}
\end{notation}

\begin{remark}\label{remark_symm_element}
	If $\Lambda=\OO G$ for some finite group $G$ (or a block thereof), then 
	the symmetrizing element $u$ may be chosen so that
	\begin{equation}
	    u_i=\frac{\chi_i(1)}{m_i\cdot|G|} \in \Q^\times
	\end{equation}
	where $\chi_i$ is the $i$-th irreducible $K$-character of $G$ (or in the block under consideration),
	and $m_i$ is the number of absolutely irreducible characters it splits up into when passing from
	$K$ to its algebraic closure $\bar K$ (see Remark \ref{remark_symm_elem}).
	In particular two of the $u_i$ are equal if (and only if) the corresponding 
	absolutely irreducible characters have equal degree. The equality of two rows 
	in the decomposition matrix is a sufficient criterion for the corresponding characters to have equal degree, and therefore for the 
	corresponding $u_i$ to be equal. Note that we potentially have two equal rows in the decomposition matrix of the principal block of  $\OO\SL_2(p^f)$ if $p$ is odd (to be precise, this happens if $f$ is even). 
\end{remark}

\begin{thm}[Unique lifting]\label{thm_unique_lift_kdelta22f}
	Let $A$ be a finite-dimensional semisimple $K$-algebra with $\dim_K Z(A) = \dim_{\bar k} Z(B_0(\bar k\Delta_2(p^f)))$. Assume $A$ is split by some
     totally ramified extension of $K$.
	Given an element $u\in Z(A)^\times$ which has $p$-valuation 
	$-f$ in every Wedderburn component of $Z(\overline{K} \otimes A)$,
	there is, up to conjugacy, at most one full $\OO$-order $\Lambda_u\subset A$ satisfying the following conditions:
	\begin{enumerate}
	\item $\Lambda_u$ is self-dual with respect to $T_u$.
	\item $k\otimes \Lambda_u$ is a split $k$-form of $B_0(\bar k\Delta_2(p^f))$ 
	\end{enumerate}
	\textbf{Addendum to the theorem (concerning the dependence on $u$):} Assume $u$ and $u'$ are two symmetrizing elements subject to the above conditions, such that $\Lambda_u$
	and $\Lambda_{u'}$ both exist. Then:
	\begin{enumerate}
	 \item If $p=2$: $\Lambda_u$ and $\Lambda_{u'}$ are conjugate.
	 \item If $p\neq 2$ and $K$ splits $A$: Let $\kappa=\frac{p^f-1}{2}$. If $\frac{u_{\kappa+1}}{u_{\kappa+2}} = \frac{u'_{\kappa+1}}{u'_{\kappa+2}}$, then $\Lambda_u$ and $\Lambda_{u'}$ are conjugate.
	 \item If $p\neq 2$ and $K$ does not split $A$: If $u_{\kappa+1}\cdot \OO^\times = u'_{\kappa+1}\cdot \OO^\times $, then $\Lambda_u$ and $\Lambda_{u'}$ are conjugate.
	\end{enumerate}
	 where $\kappa$ is the number of isomorphism classes of simple modules in $B_0(\bar k \Delta_2(p^f))$.
	
	\begin{proof}
		We assume that we are given an order $\Lambda=\Lambda_u$ satisfying the given conditions. To prove the theorem we will 
		try to conjugate $\Lambda$ into a kind of ``standard form'' depending on $u$. We let $I'$ be an ideal in $kQ$ 
		as described in Proposition \ref{prop_kform_delta2_rels} such that $k\otimes_\OO\Lambda \iso kQ/I'$ (we will assume 
		that we have fixed an isomorphism and identify the two). 
		Also, as before, we denote the idempotents in $kQ$ by $e_{2i}$ and the arrows by $s_{2i,q}$.
		We wish to treat the case where $K$ splits $A$ and the case where $K$ does not split $A$
		as well as the cases $p$ even and $p$ odd  
		(essentially) uniformly. So assume that
		\begin{equation}\label{eqn_dhuiehuie}
		  A = \left(\bigoplus_{i=1}^{\kappa} K\right) \oplus \tilde{K}^{\kappa \times \kappa} \quad\textrm{ with } \kappa=\left\{ 
		  \begin{array}{cc} \displaystyle\frac{p^f-1}{2} & \textrm{ if $p\neq 2$} \\ \\ 2^f-1 & \textrm{ if $p=2$} \end{array}\right.
		\end{equation}
		where $\tilde{K}$ is isomorphic to $K$ if $p=2$, to $K\oplus K$ if $p\neq 2$ and $A$ is $K$-split, or to a fully ramified extension of $K$ of degree two if $p\neq 2$ and $A$ is not $K$-split.
		By $\tilde{\eps}$ denote the unit element of $\tilde{K}$, construed as an idempotent in $Z(A)$. 
		For each $i$ let $\hat{e}_{2i}\in\Lambda$ be a lift of $e_{2i}\in kQ/I'$, and assume without loss that 
		$\tilde\eps \hat{e}_{2i}$ is the $i$-th diagonal idempotent in $\tilde K^{\kappa\times\kappa}$ (this 
		may certainly be achieved by conjugating $\Lambda$ by an element of $A^\times$).
		Assume furthermore that $(1-\tilde \eps )\cdot\hat e_{2i}$ has non-zero entry in the 
		$i$-th  direct summand of the decomposition (\ref{eqn_dhuiehuie}). Hence we have fixed the elements $\hat e_{2i}$ as elements of the algebra $A$ as described in (\ref{eqn_dhuiehuie}).
		 Now, using the fact that $\Lambda$ is supposed to be symmetric with respect to $T_u$, it
		follows that
		\begin{enumerate}
		\item 	If $p$ is odd and $K$ splits $A$:
		 \begin{equation}
			\hat{e}_i \Lambda \hat{e}_i = \left\langle [1,1,1], [0,p^{\frac{f}{2}},-c\cdot p^{\frac{f}{2}}], [0,0,p^f] \right\rangle_\OO \subset \OO \oplus \OO \oplus \OO
			\quad \textrm{where } c= \frac{u_{\kappa+1}}{u_{\kappa+2}}
		\end{equation}
		This follows simply from the fact that a self-dual order (with respect to $T_u$) in $\OO\oplus\OO\oplus \OO$ must have elementary divisors $1,p^{\frac{f}{2}},p^f$ (as an $\OO$ lattice in $\OO\oplus \OO \oplus \OO$) and all traces with respect to $T_u$ must be integral. Note that this
		also implies that $f$ must be even (in this situation, i. e. when $K$ splits $A$ and $p$ is odd).
		\item If $p$ is odd and $K$ does not split $A$:
		\begin{equation}
			\hat{e}_i \Lambda \hat{e}_i = \left\langle [1,1], [0,c \cdot \pi^f], [0,c^2 \cdot \pi^{2f}] \right\rangle_\OO \quad \textrm{ for some $c\in \OO[\pi]^\times$}
		\end{equation}
		 where $\pi$ is some uniformizer for the integral closure of $\OO$ in $\tilde{K}$, which is a fully ramified extension of $K$ in this case. Up to this point, we have  used two facts: First, that 
		the elementary divisors of $\OO[\pi] \otimes\hat{e}_i \Lambda \hat{e}_i$ (as a lattice in $\OO[\pi] \oplus \OO[\pi] \oplus \OO[\pi]$) must be $1,\pi^f,\pi^{2f}$, and second, that $\hat{e}_i \Lambda \hat{e}_i$ is generated by a single element as an $\OO$-order (since $e_i \cdot kQ/I \cdot e_i \iso k[T]/(T^3)$ is generated by a single element as a $k$-algebra).
		 In this case we need to put in some work to show that $\hat{e}_i \Lambda \hat{e}_i$ is uniquely determined (since different choices of $c$ may give rise to different orders). Note $T_u(\{0\}\oplus p^f \OO[\pi]) \subseteq \OO$, and 
		hence necessarily $\{0\}\oplus p^f \OO[\pi] \subset (\hat{e}_i \Lambda \hat{e}_i)^\sharp =\hat{e}_i \Lambda \hat{e}_i$. Moreover an element $[0,\tilde c \cdot \pi^f]$ lies in $\hat e_i \Lambda \hat e_i$ if and only if $T_u([0,\tilde c \cdot \pi^f])\in\OO$ (the reason being: a product of $[0,\tilde c \cdot \pi^f]$ 
		with another element of the same form, i. e. an element which has a non-zero entry only in the second component, will lie in $\{0\}\oplus p^f \OO[\pi]$ and therefore will map to something integral under $T_u$; it is thus only the product of $[0,\tilde c \cdot \pi^f]$ with $[1,1]$ for which it is not clear whether it gets mapped to something integral under $T_u$). This characterizes $\hat e_i \Lambda \hat e_i$ as
		\begin{equation}\label{eqn_dshjkh7euydvvdd}
			\hat{e}_i \Lambda \hat{e}_i = \OO\left[[0,\tilde c \cdot \pi^f] \ \bigg| \ T_u([0,\tilde c \cdot \pi^f] ) \in \OO\right]
		\end{equation}
		which is obviously uniquely determined by $u$ and the extension $\tilde K / K$.
 		\item If $p=2$ then
		\begin{equation}
		 	\hat e_i \Lambda \hat e_i = \left\langle [1,1], [0,2^f] \right\rangle_\OO
		\end{equation}
 		by the same argument as in the first point.
		\end{enumerate}
		In the above considerations we have used that each $u_i$ has $p$-valuation $-f$. In the case $p=2$ we have not used any further 
		information on $u$. In the case  $p\neq 2$ we have used the value of the quotient $u_{\kappa+1}/u_{\kappa+2}$ if $K$ splits $A$ and the class $u_{\kappa+1}\cdot \OO^\times$ if it does not (since the characterization in (\ref{eqn_dshjkh7euydvvdd}) depends only on  $u_{\kappa+1}\cdot \OO^\times$; note that $u_{\kappa+1}$ is an element of $\tilde K$ in this case while in the split case
		$u_{\kappa+1}$ and $u_{\kappa+2}$  are both elements of $K$). Since we will not make any further use of
		the symmetrizing element $u$ below, this will imply the addendum on the dependence on $u$.

		Note that in either case the $\hat{e}_i\Lambda\hat{e}_i$ are equal (when we identify the unique maximal orders containing them).
		In particular the image in $\End_K(\tilde K)$ of the action homomorphism of $\hat{e}_i \Lambda \hat{e}_i$ on 
		$\hat{e}_i \Lambda \hat{e}_j \subset \tilde{K}$ is the same as the image 
		of $\hat{e}_j\Lambda\hat{e}_j$ under the corresponding action homomorphism. Hence the submodule structure of $\hat{e}_i\Lambda\hat{e}_j$ is independent of whether it is construed as a
		 left $\hat{e}_i\Lambda\hat{e}_i$-module or a right $\hat{e}_j\Lambda\hat{e}_j$-module. Now 
		$e_i \cdot kQ/I'\cdot e_j$ is free as a $e_i \cdot kQ/I' \cdot e_i / \Soc(e_i \cdot kQ/I' \cdot e_i)$ left module (this is actually best seen by using the relations over $\bar k$ as given in
		Proposition \ref{prop_rel_delta22f} and then descending to $k$),
		and since $e_i \cdot kQ/I' \cdot e_i / \Soc(e_i \cdot kQ/I' \cdot e_i) \iso k\otimes \tilde \eps \hat e_i \Lambda \hat e_i$, this implies that $\hat e_i \Lambda \hat e_j$ is free as a left $\tilde\eps\hat e_i \Lambda \hat e_i$-module.
		This implies (when $\hat e_i A \hat e_j$ is identified with $\tilde K$ in the natural way)
		\begin{equation}\label{formula78}
			\hat e_i \Lambda \hat e_j = x_{ij} \cdot \tilde\eps \hat e_i \Lambda \hat e_i \quad \textrm{ for some $x_{ij} \in \tilde K^\times	$}
		\end{equation}
		In addition, we may and will assume that the $x_{ij}$ are integral over $\OO$. For each $i$ and $q$ we have
		\begin{equation}
			\prod_{l=0}^{p-1} e_{i + l\cdot \ul q}\cdot  kQ/I'\cdot  e_{i+(l+1)\cdot \ul q} = 0	 
		\end{equation}
		and hence
		\begin{equation}\label{eqn_lhj3228bb2623hj}
			\prod_{l=0}^{p-1} \hat e_{i + l\cdot \ul q}\cdot  \Lambda \cdot  \hat e_{i+(l+1)\cdot \ul{q}}	
			\subseteq p\cdot \hat e_i\cdot \Lambda\cdot \hat e_{i+\ul {q+1}}
		\end{equation}
		Everything from here down to (\ref{eqn_fix_elambhjdhj}) below is about showing that the inclusion in (\ref{eqn_lhj3228bb2623hj}) is in fact an equality. 
		The significance of this is that it can then be used as a formula to compute the
		$\hat e_i\cdot \Lambda\cdot \hat e_{i+\ul {q+1}}$ from the $\hat e_i\cdot \Lambda\cdot \hat e_{i+\ul {q}}$,
		showing that $\Lambda$ is determined by the $\hat e_i\cdot \Lambda\cdot \hat e_{i+\ul {0}}$. 
    
		We define a ``normalized index'' for full $\OO$-lattices $L_1\supseteq L_2$ in $\tilde K$ as follows:
		\begin{equation}
			{\rm idx} (L_1,L_2) := \frac{\length_\OO L_1/L_2}{\length_\OO L_1/p L_1}
		\end{equation}
		Note that the denominator is a constant independent of the choice of $L_1$.
		For arbitrary lattices $L_1,L_2 \subset \tilde K$ (neither of which necessarily contains the other) we define 
		${\rm idx} (L_1,L_2) := {\rm idx} (L_1+L_2,L_2) - {\rm idx} (L_1+L_2,L_1)$. 
		Now, if $L$ is any full lattice in $\tilde K$, and $x_1,x_2 \in \tilde K ^\times$, then 
		\begin{equation}
		{\rm idx}(L, x_1\cdot x_2\cdot L) = {\rm idx}(L, x_1\cdot L)
		+{\rm idx}(L, x_2\cdot L)
		\end{equation} because ${\rm idx}(L,x_i\cdot L)$ equals a constant multiple of the $p$-valuation
		of the determinant of ``multiplication with $x_i$'' construed as a $K$-vector space automorphism of 
		$\tilde{K}$.
		Now define
		\begin{equation}
			m_{i,q} := {\rm idx}  \left({\tilde\eps} \hat e_i \Lambda \hat e_i, \hat e_i \Lambda \hat e_{i+\ul q} \right)
		\end{equation}
		where  we view $\hat e_i \Lambda \hat e_{i+\ul q}$ as a subset of $\tilde K$
		as in (\ref{formula78}).
		Define furthermore 
		\begin{equation}\label{eqn_defi_aiq}
			a_{i,q} := {\rm idx} \left( \hat e_i\cdot \Lambda\cdot \hat e_{i+\ul {q+1}}, \prod_{l=0}^{p-1} \hat e_{i +l\cdot \ul q}\cdot  \Lambda \cdot  \hat e_{i+(l+1)\cdot \ul {q}}\right) = \left(\sum_{l=0}^{p-1} m_{i+l\cdot \ul q,q}\right) - m_{i,q+1}
		\end{equation}
		Clearly $a_{i,q}\geq 1$ for all $i$ and $q$.
		We have for any $q\neq r$
		\begin{equation}
	 		e_i \cdot kQ/I' \cdot e_{i+\ul q} \cdot kQ/I' \cdot e_{i+\ul q+\ul r} = e_i \cdot kQ/I' \cdot e_{i+\ul q+\ul r}
	 	\end{equation}
		and hence in particular
		\begin{equation}
			\hat e_i \Lambda \hat e_{i+\ul q} \Lambda \hat e_{i+\ul{q}+\ul{q+1}} = 
			\hat e_i \Lambda  \hat e_{i+\ul{q}+\ul{q+1}} =
			\hat e_i \Lambda \hat e_{i+\ul{q+1}} \Lambda \hat e_{i+\ul{q}+\ul{q+1}}
		\end{equation}
		which implies for all $i$ and $q$ that
		\begin{equation}\label{eqn_ext_djkljji}
			m_{i,q}+m_{i+\ul q,q+1} = m_{i,q+1}+m_{i+\ul{q+1},q}		 
		\end{equation}
		Now
		\begin{equation}
			\begin{array}{rcl}
			\displaystyle a_{i,q}-a_{i+\ul q,q} &=& \displaystyle  \left(\sum_{l=0}^{p-1} m_{i+l\cdot\ul q,q}\right)  - \left(\sum_{l=1}^{p} m_{i+l \cdot \ul q,q}\right)- m_{i,q+1} + m_{i+\ul q,q+1} \\&=& m_{i,q} - m_{i+\ul {q+1},q}- m_{i,q+1} + m_{i+\ul q,q+1} \stackrel{(\ref{eqn_ext_djkljji})}{=} 0
			\end{array} 
		\end{equation}
		Since $p$ is relatively prime to $\kappa$, this implies that $a_{i,q}=a_q$ for some
		$a_q$ independent of $i$.
		%(note that we have restricted our attention to the principal block,
		%otherwise we would have to say that $a_q$ depends only on the parity of $i$).
		Now we sum up (\ref{eqn_defi_aiq}) over all $\kappa$ values of $i$, and get
		\begin{equation}
			\sum_{i=1}^\kappa m_{2i,q+1} = p \cdot \sum_{i=1}^{\kappa} m_{2i,q} - \kappa \cdot a_q
		\end{equation}
		Plugging this formula into itself $f$ times yields (for all values of $q$)
		\begin{equation}
			\sum_{i=1}^\kappa m_{2i,q} = p^f \cdot \sum_{i=1}^{\kappa} m_{2i,q} - \kappa\sum_{i=1}^{f} p^{f-i}\cdot a_{q+i-1}			
		\end{equation}
		which implies
		\begin{equation}\label{formula_91skplks}
			\sum_{i=1}^\kappa m_{2i,q} = \frac{\kappa}{p^f-1}\cdot \sum_{i=1}^{f} p^{f-i}\cdot a_{q+i-1}	\geq \frac{\kappa}{p-1}
		\end{equation}
		with equality if and only if all $a_{q}$ are equal to $1$. Now we know (by inspecting the quiver relations) that
		\begin{equation}\label{eqn_kkdhhe77eheh}
			\begin{array}{rclc}
			\displaystyle\Jac(e_i \cdot kQ/I' \cdot e_i) &=&\displaystyle \prod_{q=0}^{f-1} \prod_{j=1}^{\frac{p-1}{2}} e_{i+\frac{1}{2}\cdot (\ul q - \ul 0) + (j-1)\cdot \ul q} \cdot kQ/I' \cdot e_{i+\frac{1}{2}\cdot (\ul q - \ul 0) + j\cdot \ul q} & \textrm{ ($p\neq 2$)} \\
	\displaystyle\Jac(e_i \cdot kQ/I' \cdot e_i) &=&\displaystyle \prod_{q=0}^{f-1} e_{i+\ul q - \ul 0} \cdot kQ/I' \cdot e_{i+\ul{q+1}-\ul 0} & \textrm{ ($p=2$)} \\
			\end{array}
		\end{equation}
		In the upper equation we used that $\frac{1}{2}(\ul q - \ul 0)=\sum_{r=0}^{q-1} \frac{p-1}{2}\ul r$.
		Now 
% 		for all $i$ we have that 
% 		\begin{equation}
% 		    \hat e_i\Lambda \hat e_i = \tilde\eps \hat e_i \Lambda \hat e_i \cap \hat e_i \Lambda \hat e_i \oplus \OO\cdot \hat e_i 
% 		\end{equation}
% 		This means that
		$\tilde\eps \hat e_i \Lambda \hat e_i \cap \hat e_i \Lambda \hat e_i$ is a pure sublattice of $\hat e_i\Lambda \hat e_i$. The
		$k$-dimension of its
		image in $e_i\cdot kQ/I'\cdot e_i$ must therefore be equal to its $\OO$-rank (which is one if $p=2$ and two otherwise), which implies that said image is equal to $\Jac(e_i\cdot kQ/I'\cdot e_i)$. Another ramification of $\tilde\eps \hat e_i \Lambda \hat e_i \cap \hat e_i \Lambda \hat e_i$ being a pure sublattice of
		$\hat e_i \Lambda \hat e_i$ is that any proper sublattice of it maps to a proper subspace of $\Jac(e_i\cdot kQ/I'\cdot e_i)$.
		Hence (\ref{eqn_kkdhhe77eheh}) implies the following:
		\begin{equation}\label{form_93jldjdkl}
			\begin{array}{rclc}
			\displaystyle\tilde\eps \hat e_i \Lambda \hat e_i \cap \hat e_i \Lambda \hat e_i &=&\displaystyle \prod_{q=0}^{f-1} \prod_{j=1}^{\frac{p-1}{2}} \hat e_{i+\frac{1}{2}\cdot (\ul q - \ul 0) + (j-1)\cdot \ul q} \Lambda \hat e_{i+\frac{1}{2}\cdot (\ul q - \ul 0) + j\cdot \ul q} & \textrm{($p\neq 2$)}\\
			\displaystyle\tilde\eps \hat e_i \Lambda \hat e_i \cap \hat e_i \Lambda \hat e_i &=&\displaystyle \prod_{q=0}^{f-1} \hat e_{i+\ul q - \ul 0} \Lambda \hat e_{i+\ul{q+1} - \ul 0} & \textrm{ ($p=2$)} \\
			\end{array}
		\end{equation}
		This, in turn, implies that the following holds for any index $i$:
		\begin{equation}
			\begin{array}{rccclc}
			\displaystyle\frac{f}{2}&=&\displaystyle{\rm idx}(\tilde\eps \hat e_i \Lambda \hat e_i , \tilde\eps \hat e_i \Lambda \hat e_i \cap \hat e_i \Lambda \hat e_i)&=&\displaystyle\sum_{q=0}^{f-1}\sum_{j=1}^{\frac{p-1}{2}}m_{i+\frac{1}{2}\cdot (\ul q - \ul 0)+(j-1)\cdot\ul q,q}&\textrm{($p\neq 2$)} \\
			\displaystyle f&=&\displaystyle{\rm idx}(\tilde\eps \hat e_i \Lambda \hat e_i , \tilde\eps \hat e_i \Lambda \hat e_i \cap \hat e_i \Lambda \hat e_i)&=&\displaystyle\sum_{q=0}^{f-1} m_{i+\ul q -\ul 0,q}&\textrm{($p=2$)} 
			\end{array}
		\end{equation}
		Summing this up over all $\kappa$ different values of $i$ yields (regardless of whether
		$p$ is even or odd)
		\begin{equation}
		 	\kappa \cdot \frac{f}{2} = \sum_{q=0}^{f-1}\frac{p-1}{2}\sum_{i=1}^{\kappa} m_{2i,q}
		\end{equation}
		Now we plug in (\ref{formula_91skplks}) to get
		\begin{equation}
			\kappa \cdot \frac{f}{2} = \frac{p-1}{2}\cdot\frac{\kappa}{p^f-1}\cdot\sum_{q=0}^{f-1} \sum_{i=1}^{f} p^{f-i}\cdot a_{q+i-1} =  \frac{p-1}{2}\cdot\frac{\kappa}{p^f-1}\cdot \frac{p^f-1}{p-1}\cdot \sum_{q=0}^{f-1} a_q
		\end{equation}
		We conclude
		\begin{equation}
			\sum_{q=0}^{f-1} a_q = f
		\end{equation}
		which implies that all $a_q$ are equal to one. This implies that the
		$\hat e_{2i} \Lambda \hat e_{2i+\ul 0}$ determine $\Lambda$ in the sense that the formula
		\begin{equation}\label{eqn_fix_elambhjdhj}
			\hat e_{2i} \Lambda \hat e_{2i+\ul{q+1}} = \frac{1}{p} \cdot \hat e_{2i} \Lambda \hat e_{2i+\ul q}\cdots \hat e_{2i+(p-1)\cdot\ul q} \Lambda \hat e_{2i+p\cdot \ul q}		 
		\end{equation}
		shows how to calculate $\hat e_{2i}\Lambda \hat e_{2i+\ul{q+1}}$ from the knowledge of
		the $\hat e_{2j}\Lambda \hat e_{2j+\ul q}$ (for all $j$).

		Now we may replace $\Lambda$ by $y^{-1}\cdot \Lambda \cdot y$, where
		\begin{equation}
			y := \left[ 1,\ldots,1, \diag\left(\prod_{j=0}^{i-1} x_{2j,2j+\ul 0} \ \bigg| \ i=1,\dots,\kappa\right) \right] \in A^\times 
		\end{equation}
		(the $x_{ij}$ were defined in (\ref{formula78}))
		and so we may assume without loss that all $x_{2i,2i+\ul 0}$ are equal to $1$, except possibly 
		$x_{2\kappa-\ul 0,2\kappa}$. 
		In other words, we have fixed all but one of the $\hat e_{2i} \Lambda \hat e_{2i+\ul 0}$. But we have
		\begin{equation}\label{eqn_hhhdh3ghjkgxjhd}
			\hat e_{2\kappa-\ul 0} \Lambda \hat e_{2\kappa} = \left\{ v \in  \hat e_{2\kappa-\ul 0} A\hat e_{2\kappa}
			\ | \  \hat e_{2\kappa} \Lambda \hat e_{2\kappa-\ul 0} \cdot v \subseteq \hat e_{2\kappa} \Lambda \hat e_{2\kappa} \right\}
		\end{equation}
		which is a consequence of the fact that $\hat e_{2\kappa-\ul 0} \Lambda \hat e_{2\kappa} $ is the dual of $\hat e_{2\kappa} \Lambda \hat e_{2\kappa-\ul 0} $ with respect to the
		bilinear pairing induced by $T_u$ (this is a general fact on self-dual orders independent of the concrete symmetrizing form $T_u$; in fact $u$ does not even show up in (\ref{eqn_hhhdh3ghjkgxjhd})).  Now in the above formula, $\hat e_{2\kappa}\Lambda \hat e_{2\kappa}$ is explicitly known, and $  \hat e_{2\kappa} \Lambda \hat e_{2\kappa-\ul 0}$ can be calculated by repeated application of (\ref{eqn_fix_elambhjdhj}) from the $\hat e_{2i} \Lambda e_{2i+\ul 0}$ with $0 \leq i < \kappa-1$ (which were fixed above by means of conjugation). This can be seen by realizing that $e_{2\kappa}\cdot kQ/I' \cdot e_{2\kappa-\ul 0}$
		can be written as a product of various $e_{2i}\cdot kQ/I' \cdot e_{2i+\ul q}$ with $0 \leq 2i < 2i+{\ul q} \leq 2(\kappa - 1 )$ and hence $\hat e_{2\kappa}\Lambda \hat e_{2\kappa-\ul 0}$ can be written as a product of various $\hat e_{2i}\Lambda \hat e_{2i+\ul q}$ with the same restriction in $i$ and $q$. But the restriction on $i$ and $q$ ensures that these $\hat e_{2i}\Lambda \hat e_{2i+\ul q}$ can be computed by means of (\ref{eqn_fix_elambhjdhj}) using only those $\hat e_{2i} \Lambda e_{2i+\ul 0}$ with $0 \leq i < \kappa-1$ .
		Hence, $\Lambda$ is determined in the sense that we have conjugated $\Lambda$ to some fixed order
		determined by the data given in the statement of the theorem. This concludes the proof.
    \end{proof}
\end{thm}

\begin{remark}\label{remark_lift_aut_sl2ppf}
	Situation as in the last theorem. Assume furthermore that the (unique) lift $\Lambda=\Lambda_u$ exists.
	Then the above proof also implies the following: If $\alpha \in \Aut_k(k\otimes\Lambda)$ is an automorphism
	of $k\otimes\Lambda$ permuting the set of idempotents $\{e_i\}_i$, then there exists an element
	$\hat \alpha \in \Aut_\OO(\Lambda)$ inducing the corresponding permutation on the set of idempotents
	$\{\hat e_i\}_i$. This follows simply from the fact that we fixed the idempotents at the beginning of the 
	proof of the Theorem and then only used conjugation by elements of $A^\times$ that commuted with all $\hat e_i$ to conjugate $\Lambda$ to any potential other lift of $k\otimes\Lambda$ (also containing the same fixed set of idempotents $\{ \hat e_i \}_i$).
\end{remark}
% 
% The following Corollary should be understood as an intermediate result, 
% the final (stronger) version of which will be Corollary \ref{corollary_unique_lift_sl2}.
% \begin{corollary}
% 	Let $\overline{\Lambda}$ be a $k$-algebra which is derived equivalent to $kQ/I$.
% 	Then for any given finite-dimensional semisimple $K$-algebra $A$ with $\dim_K Z(A) = \dim_k Z(\overline\Lambda)$, and any given element $u\in Z(A)$ which has $p$-valuation 
% 	$-f$ in every Wedderburn component of $Z(\overline{K} \otimes A)$,
% 	there is at most one derived equivalence class of  $\Lambda \in \mathfrak{L}_s(\overline\Lambda)$ with $K\otimes \Lambda \iso A$ which is self-dual with respect to $T_u$.
% \end{corollary}

\section{Transfer to $\OO{\rm SL}_2(p^f)$} 
Now we will generalize the result of Theorem \ref{thm_unique_lift_kdelta22f} to all algebras derived equivalent to a split $k$-form of $B_0(\bar k\Delta_2(p^f))$.  This will in particular include the 
two non-semisimple blocks of $k\SL_2(p^f)$. 

\begin{lemma}\label{lemma_outs_equal_out0}
	Let $k$ be algebraically closed and let $B$ be the principal block of $k\Delta_2(p^f)$.
	There is an epimorphism of algebraic groups
	\begin{equation}\label{epi_auto_grps_sl}
		\prod_{i=1}^{f} Z(B)^\times \twoheadrightarrow \Out^s_k(B)
	\end{equation}
	In particular, $\Out^s_k(B)$ is connected as an algebraic group, and hence equal to  $\Out^0_k(B)$.
\begin{proof}
	We retain the notations of the previous section, and in particular we identify $B$ with a block of $kQ/I$ (with $Q$ and $I$ as defined in Proposition \ref{prop_rel_delta22f}).
	First define a homomorphisms of algebraic groups 
	\begin{equation}\label{eqn_defipsihjkhjh33}
		\psi: \ \prod_{i=1}^{f} Z(B)^\times \rightarrow \Aut_k^s(B)
	\end{equation}
	which sends $(z_1,\ldots,z_f)$ to the automorphism given by $s_{i,q} \mapsto z_q \cdot s_{i,q}$
	(and mapping the $e_i$ to themselves).
	It is clear that those are automorphisms by checking that the images satisfy the relations given in Proposition
	\ref{prop_rel_delta22f}. We claim that the composition of $\psi$ with the natural epimorphism
	$\Aut_k^s(B) \twoheadrightarrow \Out_k^s(B)$ is surjective.	Note that $Z(B)^\times$ is an extension of $\mathbb{G}_m(k)$ by the affine plane $\Jac(Z(B))$, and hence is connected.

	We first prove the following claim, which will be used below: If $n\in \N$ is relatively prime to $p$, then
	the equation $T^n - z$ for $z\in Z(B)^\times$ has a solution in $Z(B)^\times$. This follows from the fact that
	a full set of $n$ orthogonal primitive idempotents can be lifted from $k[T]/(T^n-\overline z)$ to $Z(B)[T]/(T^n-z)$ (where $\overline z$ is the image of $z$ in $Z(B)/\Jac(Z(B))=k$). This yields a 
	decomposition of algebras $Z(B)[T]/(T^n-z) \iso A_1 \oplus \ldots \oplus A_n$. Since the $A_i$ are, in particular, $Z(B)$-modules, and $Z(B)[T]/(T^n-z)$ is free of rank $n$ as a $Z(B)$-module, we must have
	that each $A_i$ is a $Z(B)$-algebra that is free of rank one as a $Z(B)$-module. Hence each $A_i$ is canonically isomorphic (as a $k$-algebra) to $Z(B)$, and the image of $T$ in any of the  $A_i \iso Z(B)$ is a solution of $T^n-z=0$.

	Now we come to the actual proof of surjectivity of the composition of $\psi$ with the natural epimorphism
		$\Aut_k^s(B) \twoheadrightarrow \Out_k^s(B)$. Assume that $\alpha \in \Aut(B)$ is an automorphism such that $P\otimes {_\id A_\alpha} \iso P$
	for all projective indecomposables $P$. All full sets of orthogonal primitive idempotents in $B$ are conjugate (see, for instance, \cite[Introduction \S6, Exercise 14]{CurtisReinerI}), and hence we may compose $\alpha$
	with an inner automorphism of $B$ such that the resulting automorphism fixes all idempotents. We replace $\alpha$  by this new automorphism (without loss of generality). Since the canonical map $Z(B) \rightarrow e_i B e_i$ is surjective,
	and $s_{i,q}$ is a generator for the $e_i B e_i$ module $e_i Be_{i+\ul q}$, we will have $\alpha(s_{i,q})= z_{i,q}\cdot s_{i,q}$ for certain elements $z_{i,q} \in Z(B)^\times$ (and the $z_{i,q}$ determine $\alpha$).
	Now consider conjugation with elements $v$ of the form  $v=\sum_i c_i e_i$ for certain $c_i \in Z(B)^\times$: 
	\begin{equation}
		v^{-1} \cdot \alpha(s_{i,q}) \cdot v = \underbrace{\frac{c_{i+\ul q}}{c_i}\cdot z_{i,q}}_{=: \tilde z_{i,q}}\cdot s_{i,q}
	\end{equation}
	With $\tilde z_{i,q}$ defined as in the above equation we have
	\begin{equation}\label{eqn_kjdjhkdjkj_iudid}
		\prod_{i} \tilde z_{i,0} = \prod_{i} z_{i,0}
	\end{equation}
	Furthermore we can choose the $c_i$  in the definition of $v$ to assign prescribed values to all but one of the $\tilde z_{i,0}$. Choose 
	the $c_i$ so that all but possibly one of the  $\tilde z_{i,0}$ become equal to an $\kappa$-th root of the above product (where $\kappa$
	is the number of simple modules in the block, which is relatively prime to $p$). Then by the invariance of the
	product given in (\ref{eqn_kjdjhkdjkj_iudid}), all $\tilde z_{i,0}$ will be equal. Replace (without loss) $\alpha$ by the composition of $\alpha$
	with conjugation by this $v$, that is, assume that all $z_{i,0}$ are equal. We claim that this $\alpha$
	(which differs from the $\alpha$ we started with only by an inner automorphism) lies in $\Im(\psi)$ (with $\psi$ as defined in (\ref{eqn_defipsihjkhjh33})).
	To show this first notice that for $q\neq r$ the product $s_{i,q}\cdot s_{i+\ul q,r}$ is a generator for the
	$e_iBe_i$-module $e_i B e_{i+\ul q+\ul r}$, which is isomorphic to the $e_iBe_i$-module $e_i Be_{i+\ul q}$.
	Hence for any $c,\tilde c \in Z(B)^\times$ we have $c\cdot s_{i,q} = \tilde c\cdot s_{i,q}$ if and only if
	$c\cdot s_{i,q}s_{i+\ul q,r} = \tilde c\cdot s_{i,q}s_{i+\ul q,r}$. Furthermore,  in order for $\alpha$ to be an automorphism, the following relation must hold: 
	\begin{equation}
		\begin{array}{rcl}
		z_{i,q}\cdot z_{i+\ul q,q+1} \cdot s_{i,q}s_{i+\ul q,q+1} &=& 
		z_{i,q+1}\cdot z_{i+\ul{q+1},q} \cdot s_{i,q+1}s_{i+\ul{q+1},q}\\ &=& 
		z_{i,q+1}\cdot z_{i+\ul{q+1},q} \cdot s_{i,q}s_{i+\ul q,q+1}
		\end{array}
	\end{equation}
	So if we assume (as an induction hypothesis) that all $z_{i,q}$ (for some fixed value of $q$) are equal, then this implies
	that $z_{i+\ul q,q+1}\cdot s_{i,q} = z_{i,q+1}\cdot s_{i,q}$, and hence we may set
	$z_{i+\ul q,q+1} =  z_{i,q+1}$. Consequentially, all $z_{i,q+1}$ are equal. Therefore $\alpha$
	agrees with an element of $\Im(\psi)$ on the generators $s_{i,q}$. But this implies $\alpha \in \Im(\psi)$. 
\end{proof}
\end{lemma}

\begin{remark}
      By determining the kernel of the epimorphism in (\ref{epi_auto_grps_sl}) one can easily deduce that
      \begin{equation}
	 \Out_k^s(B) \iso \prod^fk[T]/(T^2)^\times \iso (\mathbb G_m^f\times \mathbb G_a^f)(k) \quad \textrm{ if $p\neq 2$} 
      \end{equation}
      and
      \begin{equation}
	 \Out_k^s(B) \iso \mathbb G_m^f(k)\quad \textrm{ if $p = 2$} 
      \end{equation}
%       It is also not too hard to see what the entire outer automorphism group has to look like:
%       \begin{equation}
% 	    \Out_k(B) \iso \left((\mathbb G_m \times \mathbb G_a)(k) \wr C_f\right)\times C_{\frac{p^f-1}{2}}\quad \textrm{ if $p\neq 2$}
%       \end{equation}
%       and
%       \begin{equation}
% 	    \Out_k(B) \iso \left( \mathbb G_m(k) \wr C_f\right)\times C_{2^f-1}\quad \textrm{ if $p\neq 2$}
%       \end{equation}
\end{remark}

\begin{lemma}\label{lemma_prop_unique_applicale}
	Let $\overline\Lambda$ be a split $k$-form of the principal block $\bar k\Delta_2(p^f)$, and assume there is a lift $\Lambda$ 
	of $\overline{\Lambda}$ subject to conditions as in Theorem \ref{thm_unique_lift_kdelta22f} (by said theorem, this lift will be unique). Then if $\alpha \in \Aut_k(\overline\Lambda)$, then there
	exists a $\beta \in \Aut_\OO(\Lambda)$ such that $\alpha\circ \overline{\beta} \in \Aut_k^s(\overline\Lambda)$
	(where $\overline\beta$ denotes the image of $\beta$ in $\Aut_k(\overline\Lambda)$).
\begin{proof}
	This follows from the fact that (since any two full sets of orthogonal primitive idempotents are conjugate) the automorphism $\alpha$ can be composed with an inner automorphism (which clearly fixes all simple modules) to get an automorphism of $\overline\Lambda$ that induces a permutation on some full set
	of orthogonal primitive idempotents in $\overline\Lambda$. Now Remark \ref{remark_lift_aut_sl2ppf} implies the
	existence of $\beta$.
\end{proof}
\end{lemma}

\begin{corollary}\label{corollary_unique_lift_sl2}
      	Let $\overline\Gamma$ be a $k$-algebra that is derived equivalent
	to a split $k$-form $\overline\Lambda$ of $B_0(\bar k\Delta_2(p^f))$. Moreover let $B$ be a finite-dimensional semisimple $K$-algebra with $\dim_K Z(B) = \dim_{\bar k} Z(B_0(\bar k\Delta_2(p^f)))$ and assume  $B$ is split by some
     totally ramified extension of $K$.
	Given an element $u\in Z(B)^\times$ which has $p$-valuation 
	$-f$ in every Wedderburn component of $Z(\overline{K} \otimes B)$,
	there is, up to isomorphism, at most one full $\OO$-order $\Gamma_u\subset B$ satisfying the following conditions:
	\begin{enumerate}
	\item $\Gamma_u$ is self-dual with respect to $T_u$.
	\item $k\otimes \Gamma_u$ is isomorphic to $\overline\Gamma$.
	\end{enumerate}
\begin{proof}
	 Recall the result of Proposition \ref{prop_unique_lift_class_ext}, which stated that if $\Lambda$ is a lift of $\overline\Lambda$ for
	 which every outer automorphism of $\overline\Lambda$ may be written as a composition of (the reduction of) an automorphism of $\Lambda$
	 and an element the $\bar k$-linear extension of which lies in $\Out^0_{\bar k}(\bar k\otimes_k \overline\Lambda)$, then $\Lambda$ corresponds to a single equivalence class of lifts in $\Lifts(\overline\Lambda)$.
	 This proposition is applicable to $\overline\Lambda$ and the unique lift $\Lambda$
	 of $\overline{\Lambda}$ subject to conditions as in Theorem \ref{thm_unique_lift_kdelta22f}, since we have verified in Lemma \ref{lemma_outs_equal_out0} and Lemma \ref{lemma_prop_unique_applicale} above
	 that the conditions of the proposition are met. Theorem \ref{thm_bijection_with_conditions} shows that the equivalence classes in $\Lifts(\overline\Lambda)$
	 subject to the conditions of Theorem \ref{thm_unique_lift_kdelta22f} (with a modified $u$, depending on the choice of the derived equivalence; see Theorem \ref{thm_symm_der_eq}) are in bijection with the equivalence classes in $\Lifts(\overline\Gamma)$ subject to
	 the conditions given in the statement of this corollary.
	 Therefore there is at most one equivalence class of lifts of $\overline\Gamma$ satisfying our assumptions. In particular there is at most one isomorphism class
	 of orders satisfying the assumptions.
\end{proof}
\end{corollary}

\begin{remark}
	Brou\'{e}'s abelian defect conjecture states the following:
	Let $k$ be an algebraically closed field, $G$ a group, $B$ a block of $kG$, $P$ a defect group of $B$, and $b$ the Brauer correspondent
	of $B$ in $kN_G(P)$. Then $b$ and $B$ are derived equivalent.

	Brou\'{e}'s conjecture has been proven (in defining characteristic) for the principal block of  $\Sl_2(q)$
	in \cite{Okuyama} (although this paper has unfortunately never been published). It has
	also been shown to hold for the unique non-principal block of maximal defect of $\Sl_2(q)$
	(which exists if $q$ is odd)
	in \cite{Yoshii}.
\end{remark}

\begin{corollary}
	Assume $k$ is algebraically closed.
	Then the generators for a basic order of $\OO\Sl_2(p^f)$ as conjectured in
	\cite{NebeSl2Char2} (for $p=2$) respectively in \cite{NebeSl2Charp} (for $p$ odd) define
	an $\OO$-order which is Morita equivalent to $\OO \SL_2(p^f)$. This is because  Corollary \ref{corollary_unique_lift_sl2} holds for the blocks of $k\SL_2(p^f)$ (due to the abelian defect conjecture), guaranteeing unique lifting.
\end{corollary}

\section{Rationality of tilting complexes}

Our goal in this section is to perform a ``Galois descent for derived equivalences'' to the degree 
up to which this is possible. This will allow us to state a unique lifting theorem for
the group ring $\F_{p^f}\SL_2(p^f)$, thus ridding us of the necessity to assume an algebraically closed coefficient field.

Concerning notation: In this section we often use field extensions $\tilde K$ and $K'$ of $K$. We will always assume
that $\tilde K$ and $K'$ are (possibly infinite) algebraic extensions of $K$ of finite ramification. We denote
by $\tilde O$ respectively $\OO'$ the corresponding discrete valuation rings and by $\tilde k$ respectively $k'$
their respective residue fields. 

\begin{defi}
	We call an $\OO$-order $\Lambda$ \emph{split} if
	the $k$-algebra $k\otimes\Lambda$ is split and the $K$-algebra $K\otimes\Lambda$ is split.
\end{defi}

\begin{lemma}
    Let $k$ be finite.
    Let $\Lambda$ be an $\OO$-order such that $K\otimes\Lambda$ is split semisimple.
    Assume that there is a field extension $\tilde K/K$ of finite degree such that $\tilde\OO\otimes\Lambda$ is split and its decomposition matrix has full row rank (that is, its rank is equal to its number of columns). Then $\Lambda$ is already split.
\begin{proof}
    Assume $S$ is a simple $\Lambda$-module that is not absolutely irreducible. Since there are no non-commutative finite-dimensional 
    division algebras over $k$, $\End(S)$ is commutative and hence $\End(\tilde k\otimes S) \iso \tilde k\otimes\End(S)$ is a direct sum
    of copies of $\tilde k$.
    Therefore $\tilde k\otimes S$ is a direct sum of non-isomorphic simple $\tilde\OO\otimes\Lambda$-modules $\tilde S_1,\ldots, \tilde S_l$ (for some $l>1$). Each simple $\tilde K\otimes \Lambda$-module is of the form $\tilde K\otimes V$ for some simple
    $K\otimes\Lambda$-module $V$. Let $L$ be a $\Lambda$-lattice in $V$. Then $\tilde\OO\otimes L$ is a $\tilde\OO\otimes\Lambda$-lattice
    in $\tilde K\otimes V$, and the multiplicities of $\tilde S_1,\ldots \tilde S_l$ in $\tilde k\otimes L$ are all equal to the
    multiplicity of $S$ in $k\otimes L$. Therefore, the columns in the decomposition matrix of $\tilde\OO\otimes\Lambda$ associated to the simple modules
    $\tilde S_1,\ldots,\tilde S_l$ are all equal, in contradiction to the assumption that the decomposition matrix of $\tilde\OO\otimes\Lambda$ has full row rank. Therefore all simple $\Lambda$-modules are absolutely simple, that is, $\Lambda$ is split.
\end{proof}
\end{lemma}

\begin{lemma}\label{lemma_unramified_ext_split}
    Assume that $\tilde K$ is totally ramified over $K$. If $\Lambda$ is an $\OO$-order such that $\tilde k\otimes \Lambda$ is split,
    then $k\otimes\Lambda$ is split.
    
    In particular, under the assumption that $k$ is finite,  $\tilde K\otimes\Lambda$ is split semisimple and the decomposition matrix of $\Lambda$ over a splitting system has full row rank, $k\otimes\Lambda$ will be split.
\begin{proof}
    This is clear since $\tilde k=k$.
\end{proof}
\end{lemma}

\begin{remark}
	We should note that
	\begin{enumerate}
	 \item Full row rank of the decomposition matrix is implied if the Cartan matrix of an algebra 
		is non-degenerate (which is a known fact in the case of group rings).
	 \item Up to signs, the determinant (and therefore non-degeneracy) of the Cartan matrix is preserved under
		derived equivalences (even under stable equivalences of Morita type).
	\end{enumerate}
\end{remark}

\begin{defi}
    Let $A$ be a ring.
    We say a tilting complex $T\in\mathcal C^b(\projC_A)$ is \emph{determined by its terms}, if any tilting complex 
    $T'\in\mathcal C^b(\projC_A)$ with $T^i\iso T'^i$ for all $i\in\Z$ is isomorphic to $T$ in $\mathcal K^b(\projC_A)$.
\end{defi}

\begin{remark}
	By \cite[Corollary 8]{JensenXuZDegenerations} two-term tilting complexes defined over algebras over a field are determined by their terms. By unique lifting
	of tilting complexes (see \cite{RickardLiftTilting}), the same is true for two-term tilting complexes defined over orders over complete discrete valuation rings.
\end{remark}

\begin{defi}
    Let $\tilde \Lambda$ be an $\tilde\OO$-order. We call an $\OO$-order $\Lambda\subseteq \tilde\Lambda$ an \emph{$\OO$-form of
    $\tilde\Lambda$} if $\rank_\OO \Lambda = \rank_{\tilde\OO}\tilde\Lambda$ and $\tilde\OO\cdot \Lambda =\tilde\Lambda$.
    We define a $k$-form of a finite-dimensional $\tilde k$-algebra is the analogous way.
\end{defi}

\begin{lemma}\label{lemma_structure_complex_ring_ext}
    Let $\Lambda$ be an $\OO$-order and let $\tilde K$ be an unramified finite extension of $K$. Furthermore, let $\tilde C\in\mathcal C^b(\modC_{\tilde\OO\otimes\Lambda})$ be a complex of $\tilde \OO\otimes\Lambda$-modules and let
    $C$ be the restriction of $\tilde C$ to $\Lambda$. Then, in the category $\mathcal C^b(\modC_{\tilde\OO\otimes\Lambda})$,
    \begin{equation}
	   \tilde\OO\otimes C \iso \bigoplus_{i=1}^{[\tilde K:K]} \tilde C^{\alpha_i} 
    \end{equation}
    for certain $\alpha_i \in \Aut_\OO(\tilde\OO)$. Here, for an $\alpha\in\Aut_\OO(\tilde\OO)$, $\tilde C^\alpha$ denotes the 
    complex of $\tilde\OO\otimes \Lambda$-module the terms of which are (as sets) equal to the terms of $\tilde C$, with differential
    equal to that of $\tilde C$, but with the following twisted action of $\tilde\OO\otimes\Lambda$ on the terms:
    \begin{equation}
	   \tilde C^i\times \tilde\OO\otimes\Lambda \longrightarrow \tilde C^i: (m,a\otimes b) \mapsto m\cdot \alpha(a)\otimes b
    \end{equation}
    We claim furthermore that at least one of the $\alpha_i$ may be chosen to be the identity automorphism of $\tilde\OO$.
\begin{proof}
    First note that $\tilde\OO \otimes_\OO\tilde\OO \iso \bigoplus^{[\tilde K : K]} \tilde\OO$, since $\tilde K$ is unramified over
    $K$. For $i\in\{1,\ldots,[\tilde K:K]\}$ denote by $\eps_i$ the epimorphism from $\tilde\OO\otimes_\OO\tilde\OO$ to 
    $\tilde \OO$ given by projection to the $i$-th component of $\bigoplus^{[\tilde K : K]}\tilde\OO$ (of course, the ordering of the $\eps_i$ is not canonical). By abuse of notation, we also denote by $\eps_i$ the unique primitive idempotent in $\tilde\OO\otimes_\OO\tilde\OO$
    that gets mapped to $1$ under the projection $\eps_i$. Now we consider the complex of  $\tilde\OO\otimes_\OO\tilde\OO\otimes_\OO\Lambda$-modules $\tilde\OO\otimes_\OO \tilde C$.
    We can decompose this complex as follows:
    \begin{equation}
	   \tilde\OO\otimes_\OO \tilde C = \bigoplus_{i=1}^{[\tilde K:K]} \tilde\OO\otimes_\OO \tilde C \cdot (\eps_i\otimes 1_\Lambda)
    \end{equation}
    Now consider the embedding 
    \begin{equation}
	   \eta: \tilde\OO \hookrightarrow \tilde\OO\otimes_\OO\tilde\OO: a \mapsto a\otimes 1
    \end{equation}
    If we turn
    $\tilde\OO\otimes_\OO \tilde C$ into a complex of $\tilde\OO\otimes\Lambda$-modules via the embedding $\eta\otimes \id_\Lambda$ we get,
    by definition, $\tilde\OO\otimes_\OO C$.
    If we turn $\tilde\OO\otimes_\OO \tilde C \cdot (\eps_i\otimes 1_\Lambda)$ into a complex of $\tilde\OO\otimes\Lambda$-modules via the embedding $\eta\otimes \id_\Lambda$ we get $\tilde C^{\eps_i\circ\eta}$. So the our first claim follows (with $\alpha_i := \eps_i\circ\eta$).
    As for the claim that one of the $\alpha_i$ may be chosen equal to the identity, just note that there is an epimorphism
    $\tilde\OO\otimes_\OO\tilde\OO \twoheadrightarrow \tilde\OO: a\otimes b \mapsto a\cdot b$. Since the $\eps_i$ are in fact 
    all epimorphisms from $\tilde\OO\otimes_\OO\tilde\OO$ to $\tilde\OO$,  this epimorphism needs to be equal to some $\eps_i$.
    But then $\alpha_i=\id$.
\end{proof}
\end{lemma}

\begin{prop}[Reduction to finite field extensions]\label{prop_finite_tilting_red}
	Let $\Lambda$ and $\Gamma$ be two $\OO$-orders such that
	$\tilde\OO\otimes\Lambda$ and $\tilde\OO\otimes\Gamma$ are derived equivalent,
	and let $\tilde T$ be a tilting complex over $\tilde\OO\otimes\Lambda$ with endomorphism ring
	$\tilde\OO\otimes\Gamma$. Then there exists a finite extension $K'$ of $K$ which is contained in
	$\tilde K$ such that $\OO'\otimes \Lambda$ is derived equivalent to an $\OO'$-form $\Gamma'$ of $\tilde\OO\otimes\Gamma$, and there
	is a tilting complex $T'$ over $\OO'\otimes\Lambda$ with endomorphism ring $\Gamma'$ such that
	$\tilde\OO\otimes_{\OO'}T'\iso \tilde T$ in $\mathcal K^b(\projC_{\tilde\OO\otimes\Lambda})$.
\begin{proof}
%	As a complex, $\tilde T$ is determined uniquely by a projective $\Z$-graded module $\tilde P$
%	with an endomorphism $\tilde d$ of degree $-1$. 
%	Getting even more elementary, $\tilde P$ is given by a map $n: \Z \rightarrow \Z_{\geq 0}$ such that
%	$n^{-1}(\Z_{>0})$ is finite and an $\OO$-algebra homomorphism 
%	\begin{equation}
%		\varphi: \ \Lambda \longrightarrow \bigoplus_{i\in\Z} \tilde\OO^{n(i)\times n(i)}
%	\end{equation}
%	Also, from this point of view,
%	\begin{equation}
%		\tilde d \in \bigoplus_{i\in \Z} \tilde\OO^{n(i+1)\times n(i)} 
%	\end{equation}
%	Now it is clear that since $\varphi(\Lambda)$ has finite $\OO$-rank, there is a finite field extension
%	$K'$ of $K$ such that all entries of $\varphi(\lambda_i)$ for an $\OO$-basis $(\lambda_i)_i$ of $\Lambda$
%	and all entries of $\tilde d$ lie in $K'$. Thus we get a complex $T'\in\mathcal C^b(\modC_{\OO'\otimes\Lambda})$
%	with $\tilde\OO\otimes_{\OO'} T' \iso \tilde T$. By Lemma \ref{lemma_ring_ext_tilting_is_tilting}
%	it follows that $T'$ must be a tilting complex. 
%	
	There is some invertible complex $\tilde X\in \mathcal D^b((\tilde \OO\otimes \Lambda)^{\op}\otimes_{\tilde{\OO}} (\tilde\OO\otimes \Gamma))$ with inverse
	 $\tilde Y\in \mathcal D^b((\tilde \OO\otimes \Gamma)^{\op}\otimes_{\tilde \OO} (\tilde\OO\otimes\Lambda))$ such that the restriction of $\tilde Y$ to 
	 $\tilde \OO \otimes \Lambda$ is isomorphic to $\tilde T$ in $\mathcal D^b(\tilde{\OO}\otimes \Lambda)$.
	 We can find a finite extension $K'$ of $K$ (contained in $\tilde K$) such that there are bounded
	 complexes $X'$ and $Y'$ such that $\tilde{\OO} \otimes_{\OO'}X' \iso \tilde{X}$ and $\tilde{\OO} \otimes_{\OO'}Y' \iso \tilde{Y}$.
	 This is simply because $\tilde X$ and $\tilde Y$ can be represented by bounded complexes of finitely generated modules, and 
	 so $K'$ needs only be big enough for all terms of these complexes to be defined over $\OO'$ and for the differentials (which are made up of finitely many homomorphisms) to be defined.
	 Looking at the construction of the derived tensor product, it is clear that 
	 \begin{equation}\label{eqn_djkjhdiu33jkekjbhhkjbh3e}
	 \tilde \OO \otimes_{\OO'}^{\mathbb L} (X'\otimes_{\OO'\otimes\Gamma}^{\mathbb L} Y')\iso \tilde X\otimes_{\tilde \OO\otimes\Gamma}^{\mathbb L} \tilde Y
	 \textrm{ and } \tilde \OO \otimes_{\OO'}^{\mathbb L} (Y'\otimes_{\OO'\otimes\Lambda}^{\mathbb L} X')\iso \tilde Y\otimes_{\tilde \OO\otimes\Lambda}^{\mathbb L} \tilde X
	 \end{equation}
	 But the right hand terms in (\ref{eqn_djkjhdiu33jkekjbhhkjbh3e}) have homology concentrated in degree zero. This means that
	 $X'\otimes_{\OO'\otimes\Gamma}^{\mathbb L} Y'$ and $Y'\otimes_{\OO'\otimes\Lambda}^{\mathbb L} X'$ are isomorphic to stalk complexes 
	 in $\mathcal D^-((\OO'\otimes \Lambda)^{\op} \otimes_{\OO'} (\OO'\otimes \Lambda))$ respectively 
	 $\mathcal D^-((\OO'\otimes \Gamma)^{\op} \otimes_{\OO'} (\OO'\otimes \Gamma))$. Since tensoring with $\tilde \OO$ renders
	 them isomorphic to $0\rightarrow \tilde{\OO}\otimes\Lambda \rightarrow 0$ respectively $0\rightarrow \tilde{\OO}\otimes\Gamma \rightarrow 0$ it follows from the Noether-Deuring theorem for modules that they are isomorphic to $0\rightarrow \OO'\otimes\Lambda\rightarrow 0$ respectively $0\rightarrow \OO'\otimes\Gamma\rightarrow 0$. Therefore $X'$ and $Y'$ are invertible, and thus the restriction of $Y'$ to $\OO'\otimes\Lambda$ is a tilting complex $T'$ with $\tilde \OO \otimes_{\OO'} T' \iso \tilde T$.
	 
	By \cite[Theorem 2.1]{RickardDerEqDerFun} it follows that
	the endomorphism ring of $T'$ in $\mathcal D^b(\OO'\otimes\Lambda)$ is an $\OO'$-form of $\tilde\OO\otimes\Lambda$.
\end{proof}
\end{prop}

\begin{remark}
	We should mention the following (trivial) addendum to the above proposition:
	If $\tilde\OO$ splits $\Lambda$ and/or $\Gamma$, we may choose an $\OO'$ which splits
	$\Lambda$ and/or $\Gamma$. Similarly, if $\tilde k$ splits $k\otimes\Lambda$ and/or 
	$k\otimes\Gamma$, we may choose an $\OO'$ such that $k'$ (the residue field of $\OO'$)  
	splits $k\otimes\Lambda$ and/or $k\otimes\Gamma$.
\end{remark}

\begin{lemma}\label{lemma_ring_ext_tilting_is_tilting}
    Let $\Lambda$ be an $\OO$-order and let $T\in\mathcal C^b(\modC_\Lambda)$ be a complex with
    differential $d: T \longrightarrow T[-1]$. If $\tilde\OO\otimes T$ is a tilting complex
    for $\tilde\OO\otimes\Lambda$ (in particular $\tilde\OO\otimes T \in \mathcal C^b(\projC_{\tilde\OO\otimes\Lambda})$), then $T$ is a tilting complex for $\Lambda$. 
\begin{proof}
	First note that by Proposition \ref{prop_finite_tilting_red} we may assume that $\tilde K/K$ is a field extension of finite degree.
    If $M$ is a (finitely-generated) $\Lambda$-module such that $\tilde\OO\otimes M$ is a projective 
   $\tilde\OO\otimes\Lambda$-module, $M$ must itself be projective. This follows easily from the fact that
    $\tilde\OO\otimes M$ is projective if and only if it is a direct summand of some free module, and 
   so the restriction of $\tilde\OO\otimes M$, which is just a direct sum of copies of $M$, is a summand of a restriction
   of a free module, which is again a free module. This shows that $\tilde\OO\otimes T \in \mathcal C^b(\projC_{\tilde\OO\otimes\Lambda})$ implies $T\in\mathcal C^b(\projC_\Lambda)$.

    Now we show $\Hom_{\mathcal D^b(\Lambda)}(T,T[i])=0$ for all $i\neq 0$. So let $\varphi \in \Hom_{\mathcal C^b(\projC_\Lambda)}(T, T[i])$. Then there is a homotopy $h: \tilde\OO\otimes T \longrightarrow \tilde\OO\otimes T[i+1]$ such that
    $1_{\tilde\OO}\otimes \varphi = h\circ(1_{\tilde\OO}\otimes d)+(1_{\tilde\OO}\otimes d)\circ h$.
    Since for arbitrary $\Lambda$-modules $M$ and $N$ we have $\Hom_{\tilde\OO\otimes\Lambda}(\tilde\OO\otimes M, \tilde\OO\otimes N) \iso \tilde\OO\otimes \Hom_\Lambda(M, N)$, we can write
    \begin{equation}
	   h = \sum_{j=1}^{[\tilde K:K]} b_j\otimes h_j \quad\textrm{for certain $h_j: T \longrightarrow T[i+1]$}
    \end{equation}
    where $(b_1,\ldots,b_{[\tilde K:K]})$ is an $\OO$-basis of $\tilde\OO$ and, without loss, $b_1=1_{\tilde\OO}$.
    Hence
    \begin{equation}
	   b_1\otimes \varphi = \sum_{j=1}^{[\tilde K: K]} b_j \otimes (h_j\circ d + d\circ h_j)
    \end{equation}
    This implies
    \begin{equation}
	   \varphi = h_1 \circ d + d \circ h_1
    \end{equation}
    and therefore $\varphi$ is homotopic to the zero map.

    Now we show that $T$ generates $\mathcal K^b(\projC_\Lambda)$. To see this we look at the functor 
    \begin{equation}\Res: \mathcal K^-(\projC_{\tilde\OO\otimes\Lambda}) \longrightarrow \mathcal K^-(\projC_\Lambda)
    \end{equation}
    which, by definition, simply restricts the terms of the complexes from $\tilde\OO\otimes\Lambda$-modules to $\Lambda$-modules. Since this is an exact functor, and $\Res (\tilde \OO \otimes T)$ is just 
    a direct sum of copies of $T$, $\add (T) \supseteq\Res(\add(\tilde\OO\otimes T))$. But $0\rightarrow \tilde\OO\otimes\Lambda \rightarrow 0$ lies in $\add(\tilde\OO\otimes T)$, and therefore $0\rightarrow \Lambda \rightarrow 0$ 
    lies in $\add(T)$ (since $\Res(0\rightarrow \tilde\OO\otimes\Lambda\rightarrow 0)$ is isomorphic to a direct sum of copies of $0\rightarrow \Lambda \rightarrow 0$).
\end{proof}
\end{lemma}

\begin{thm}\label{thm_descent_derived}
    Assume $k$ is finite and $\tilde K$ is unramified over $K$.
    Let $\tilde \Lambda$ be an $\tilde \OO$-order such that $\tilde k\otimes\tilde\Lambda$ is split, $\tilde K\otimes \tilde\Lambda$ is semisimple
    and the decomposition matrix of $\tilde\Lambda$ over a splitting system has full row rank. Let $\tilde T\in \mathcal C^b(\projC_{\tilde\Lambda})$ be a tilting complex that is determined by its terms. Set 
    \begin{equation}\tilde\Gamma := \End_{\mathcal D^b(\tilde\Lambda)}(\tilde T)\end{equation}
    If $\Lambda$ is an $\OO$-form of $\tilde\Lambda$ such that $k\otimes \Lambda$ is split and 
    there is a totally ramified extension of $K$ that splits $K\otimes\Lambda$, then there is an $\OO$-form $\Gamma$ of $\tilde\Gamma$ with the same properties that is derived equivalent to $\Lambda$.
\begin{proof}
    Let $T$ be the restriction of $\tilde T$ to $\mathcal C^b(\projC_\Lambda)$.
    By Lemma \ref{lemma_structure_complex_ring_ext} the complex $\tilde\OO\otimes T$ is isomorphic to
    a direct sum of complexes of the form $\tilde T^\alpha$ for certain $\alpha \in \Aut_\OO(\tilde\OO)$.
    Now note that since $k\otimes \Lambda$ is split, the projective indecomposable $\tilde\Lambda$-modules $\tilde P$ are of the form
    $\tilde\OO\otimes P$ for projective indecomposable $\Lambda$-modules $P$. Therefore they are isomorphic to their Galois twists. In particular, the terms of $\tilde T^\alpha$ and $\tilde T$ are isomorphic for all $\alpha\in\Aut_\OO(\tilde\OO)$. Since $\tilde T$ is by assumption determined by its terms, we must have
   $\tilde T^\alpha \iso \tilde T$ for all $\alpha\in\Aut_\OO(\tilde\OO)$. This shows that
   $\tilde\OO\otimes T$ is a tilting complex, and therefore so is $T$ (by Lemma \ref{lemma_ring_ext_tilting_is_tilting}).  It is clear by 
   \cite[Theorem 2.1]{RickardDerEqDerFun} (or by using linear algebra) that the endomorphism ring of $T$ is an $\OO$-form of the endomorphism ring of
	$\tilde\OO\otimes T$, and of course it is derived equivalent to $\Lambda$. We have 
	\begin{equation}\tilde \OO \otimes \End_{\mathcal D^b(\Lambda)}(T) \iso \End_{\mathcal D^b(\tilde\Lambda)}(\tilde\OO\otimes T)\iso \tilde\Gamma^{[\tilde K:K]\times[\tilde K:K]}\end{equation} that is,  $\End_{\mathcal D^b(\Lambda)}(T)$ is on $\OO$-form
	of $\tilde\Gamma^{[\tilde K:K]\times[\tilde K:K]}$. This will yield on $\OO$-form of $\tilde\Gamma$ with
	the desired properties (simply by applying a Morita equivalence) once we see that $k\otimes \End_{\mathcal D^b(\Lambda)}(T)$ is split.
	Let $K'$ be a totally ramified extension of $K$ such that $K'\otimes\Lambda$ is split. Since
	$K'\otimes\End_{\mathcal D^b(\Lambda)}(T)\iso \End_{\mathcal D^b(K'\otimes\Lambda)}(K'\otimes T)$ is
	Morita equivalent to $K'\otimes\Lambda$, it follows by Lemma \ref{lemma_unramified_ext_split} that 
	$k\otimes\End_{\mathcal D^b(\Lambda)}(T)$ is split.
\end{proof}
\end{thm}

\begin{corollary}\label{corollary_descent_multiple_steps}
	The assertion of the preceding Theorem remains true if $\tilde\Lambda$ and $\tilde\Gamma$ are linked
	by a series of derived equivalences which all are afforded by tilting complexes that are determined by their terms.
\begin{proof}
	This follows by iterated application of the preceding theorem.
\end{proof}
\end{corollary}

\begin{corollary}\label{corr_sl2_der_eq_split_kform}
	Let $\OO$ be the $p$-adic completion of the maximal unramified extension of $\Q_p$.
	The blocks of defect $C_p^f$ of the group ring $\Z_p[\zeta_{p^f-1}]\SL_2(p^f)$ are derived equivalent to
	a $\Z_p[\zeta_{p^f-1}]$-form (split over $\F_{p^f}$) of their respective Brauer correspondent in $\OO\Delta_2(p^f)$ with $\Q_p[\zeta_{p^f-1}]$-span
	isomorphic to the $\Q_p[\zeta_{p^f-1}]$-span of the corresponding block of $\Z_p[\zeta_{p^f-1}]\Delta_2(p^f)$.
\begin{proof}
	The respective blocks of $k\SL_2(p^f)$ and $k\Delta_2(p^f)$ are 
	linked by a series of two-term complexes (see \cite{Okuyama} respectively \cite{Yoshii}). Hence the first claim follows 
	from Theorem \ref{thm_descent_derived} and Corollary \ref{corollary_descent_multiple_steps}. The assertion concerning the $\Q_p[\zeta_{p^f-1}]$-spans
	follows from the fact that the $\Q_p[\zeta_{p^f-1}]$-spans of the blocks 
	of $\Z_p[\zeta_{p^f-1}]\SL_2(p^f)$ and $\Z_p[\zeta_{p^f-1}]\Delta_2(p^f)$ which are Brauer correspondents are Morita equivalent. 
\end{proof}
\end{corollary}

\begin{corollary}\label{corollary_unique_lifting}
    Assume $k \supseteq \F_{p^f}$ and $B$ is a block of $k\SL_2(p^f)$ of maximal defect.
    Let $A$ be a finite-dimensional semisimple $K$-algebra with $\dim_K Z(A) = \dim_{k} Z(B)$. Assume $A$ is split by some
     totally ramified extension of $K$.
	Given an element $u\in Z(A)^\times$ which has $p$-valuation 
	$-f$ in every Wedderburn component of $Z(\bar{K} \otimes A)$,
	there is, up to conjugacy, at most one full $\OO$-order $\Lambda_u\subset A$ satisfying the following conditions:
	\begin{enumerate}
	\item $\Lambda_u$ is self-dual with respect to $T_u$.
	\item $k\otimes \Lambda_u$ is isomorphic to $B$. 
	\end{enumerate}
\begin{proof}
    By Corollary \ref{corr_sl2_der_eq_split_kform} the block $B$ is derived equivalent to a split $k$-form $\overline\Gamma$ of $B_0(\bar k \Delta_2(p^f))$.
    Thus the assertion follows directly from Corollary \ref{corollary_unique_lift_sl2}. 
\end{proof}
\end{corollary}

\begin{corollary}
    	The generators for a basic order of $\Z_p[\zeta_{p^f-1}]\Sl_2(p^f)$ as conjectured in
	\cite{NebeSl2Char2} (for $p=2$) respectively in \cite{NebeSl2Charp} (for $p$ odd) define
	a $\Z_p[\zeta_{p^f-1}]$-order which is Morita-equivalent to $\Z_p[\zeta_{p^f-1}] \SL_2(p^f)$.
\end{corollary}

As a corollary we can also prove that a discrete valuation ring version of the abelian defect conjecture holds for
$\Z_p[\zeta_{p^f-1}]\SL_2(p^f)$.
\begin{corollary}
    The non-semisimple blocks of $\Z_p[\zeta_{p^f-1}]\SL_2(p^f)$ are derived equivalent to their
    Brauer-correspondents in $\Z_p[\zeta_{p^f-1}]\Delta_2(p^f)$.
\begin{proof}
    As we have already seen, any non-semisimple block $\Gamma$ of $\Z_p[\zeta_{p^f-1}]\SL_2(p^f)$ is derived equivalent to the unique lift $\Lambda_u \subset \Q_p[\zeta_{p^f-1}]\otimes B_0(\Z_p[\zeta_{p^f-1}]\Delta_2(p^f))=: A$ of 
    a split $\F_{p^f}$-form of $B_0(\bar \F_p\Delta_2(p^f))$ with respect to some $u\in Z(A)$
    satisfying the conditions of Theorem \ref{thm_unique_lift_kdelta22f} (this is just putting 
    Corollary \ref{corr_sl2_der_eq_split_kform} and Theorem \ref{thm_unique_lift_kdelta22f} together).
    The addendum to Theorem \ref{thm_unique_lift_kdelta22f} tells us that if $p=2$, then $\Lambda_u\iso B_0(\Z_2[\zeta_{2^f-1}]\Delta_2(2^f))$ which implies the assertion of this corollary. If $p \neq 2$ and $\Q_p[\zeta_{p^f-1}]$ does not split $\Sl_2(p^f)$, then the addendum tells us  (using the same notational conventions as in Theorem \ref{thm_unique_lift_kdelta22f}, including Notation \ref{notation_cpi_ordering_sl2}; these will be used throughout this proof) that $\Lambda_u$ depends only on
    $u_{\kappa+1}\cdot \OO^\times$, which we may assume to be equal to $p^{-f}\cdot \OO^\times$ by virtue of 
    $u_{\kappa+1}$ being rational. So again, $\Lambda_u\iso B_0(\Z_p[\zeta_{p^f-1}]\Delta_2(p^f))$ follows and we are done.
    Now if $p$ is odd and $\Q_p[\zeta_{p^f-1}]$ does split $\Sl_2(p^f)$, then $\Lambda_u$ depends only on the quotient $u_{\kappa+1}/u_{\kappa+2}$. Assume for the rest of the proof that we are in this case. 
    We also fix some tilting complex $T$ in $\mathcal K^b(\projC_{\Lambda_u})$ 
    with endomorphism ring $\Gamma$.
    Furthermore let $V_{\kappa+1}$ and $V_{\kappa+2}$ be the $(\kappa+1)$-st and 
    $(\kappa+2)$-nd simple $\bar \Q_p \otimes A$-module.
    Note that the symmetrizing element $u$ for $\Lambda_u$ arises from the symmetrizing element $u'$  we use for 
    $\Gamma$ by flipping signs in certain Wedderburn components. As mentioned in Remark \ref{remark_symm_element}, $u'$ 
    may be chosen so that $u'_{\kappa+1}=u'_{\kappa+2}$, since the corresponding 
    rows in the decomposition matrix are equal (we do not make use of any particular knowledge of the decomposition matrix of $\SL_2(p^f)$ to establish this; the fact that the $(\kappa+1)$-st and $(\kappa+2)$-nd row of the decomposition matrix of $\Delta_2(p^f)$ over a splitting system are equal implies that the corresponding rows
    in the decomposition matrix of a derived equivalent order will also be equal). The sign of $u'_{\kappa+1}$ respectively $u'_{\kappa+2}$ is flipped upon passage to $\Lambda_u$ depending on
    the sign of $[V_{\kappa+1}]$  respectively $[V_{\kappa+2}]$ as a coefficient of
    \begin{equation}
	\sum_{i} (-1)^i\cdot [\bar \Q_p\otimes_{\Z_p[\zeta_{p^f-1}]} T^i] \in K_0(\modC_{\bar \Q_p\Delta_2(p^f)})
    \end{equation}
    These signs are equal, since all of the $T^i$ are projective modules and therefore 
    $V_{\kappa+1}$ and $V_{\kappa+2}$ occur in their $\bar \Q_p$-span with the same multiplicities (again since the corresponding rows in the decomposition matrix are equal). We conclude that $u_{\kappa+1}=u_{\kappa+2}$, and
    therefore $\Lambda_u \iso B_0(\Z_p[\zeta_{p^f-1}]\Delta_2(p^f))$, which is what we wanted to prove.
\end{proof}
\end{corollary}

\section*{Acknowledgements} 
During the write-up of this paper I was supported by the DFG (German Research Foundation) SPP 1388, and later the FWO Vlaanderen  (Research Foundation Flanders).
This research formed part of my PhD thesis.

\addcontentsline{toc}{section}{References}
\renewcommand\bibname{References}
\bibliographystyle{alpha}
\bibliography{refs}
\end{document}